\documentclass[12pt]{article}

\usepackage{CJK}
\usepackage{amsfonts,amssymb,amsmath,mathrsfs,multirow,booktabs}

\usepackage{color,latexsym,amsfonts}
 \newcommand{\qed}{\hfill\rule{2mm}{3mm}\vspace{4mm}}

 \newtheorem{theorem}{Theorem}[section]
 \newtheorem{lemma}[theorem]{Lemma}
 \newtheorem{corollary}[theorem]{Corollary}
 \newtheorem{proposition}[theorem]{Proposition}
 \newtheorem{example}[theorem]{Example}
 \newtheorem{Definition}[theorem]{Definition}
 \newtheorem{remark}[theorem]{Remark}
 \newtheorem{condition}[theorem]{Condition}

 \def\blemma{\begin{lemma}\sl{}\def\elemma{\end{lemma}}}
 \def\btheorem{\begin{theorem}\sl{}\def\etheorem{\end{theorem}}}
 \def\bcorollary{\begin{corollary}\sl{}\def\ecorollary{\end{corollary}}}
 \def\bdefinition{\begin{Definition}\sl{}\def\edefinition{\end{Definition}}}
 \def\bproposition{\begin{proposition}\sl{}\def\eproposition{\end{proposition}}}
 
 \def\bcondition{\begin{condition}\sl{}\def\econdition{\end{condition}}}

 \def\beqlb{\begin{eqnarray}}\def\eeqlb{\end{eqnarray}}
 \def\beqnn{\begin{eqnarray*}}\def\eeqnn{\end{eqnarray*}}

 \def\mbb{\mathbb}\def\mbf{\mathbf}

 \def\ar{&\!\!}

 \def\eqref#1{{\rm(\ref{#1})}}

 \def\proof{\noindent{\it
 Proof.~}}\def\qed{\hfill$\Box$\medskip}

\def\e{{\mbox{\rm e}}}
 \def\<{\langle}\def\>{\rangle}

 \def\mbb{\mathbb}
  \def\mbf{\mathbf}
\newcommand{\dd}{\mathrm{d}}

\newfam\msbmfam\font\tenmsbm=msbm10\textfont
\msbmfam=\tenmsbm\font\sevenmsbm=msbm7
\scriptfont\msbmfam=\sevenmsbm

\def\de{{\delta}}

\def\({\left(}\def\){\right)}

\begin{document}

\noindent{(Draft: 2020/06/01)}
\bigskip\bigskip

\centerline{\Large\bf SPDEs with non-Lipschitz coefficients
}

\smallskip\smallskip

\centerline{\Large\bf and nonhomogeous boundary conditions}

\bigskip

\centerline{
Jie Xiong\footnote{Department of Mathematics and SUSTech International center for Mathematics, Southern University of Science \& Technology, Shenzhen, China. Southern University of Science and Technology Start up fund Y01286120 and NSFC (Nos.~61873325 and 11831010).  Email: xiongj@sustech.edu.cn},
Xu Yang\footnote{School of Mathematics and Information Science,
North Minzu University, Yinchuan, China. Supported by
NSFC (No.~11771018) and Major research project for North Minzu University
(No. ZDZX201902). Email: xuyang@mail.bnu.edu.cn. Corresponding author.}}

\bigskip\bigskip

{\narrower{\narrower

\noindent{\bf Abstract.}
In this paper we establish the strong existence, pathwise uniqueness
and a comparison theorem to a stochastic partial differential equation
driven by Gaussian colored noise with non-Lipschitz
drift, H\"older continuous diffusion coefficients and
the spatial domain in finite interval, $[0,1]$, and with
Dirichlet, Neumann or mixed nonhomogeneous  random conditions imposed on the endpoints. The H\"older continuity of the solution both in time
 and in space variables is also studied.

\bigskip

\textit{Mathematics Subject Classifications (2010)}:
60H15; 60J50.

\bigskip

\textit{Key words and phrases}: Stochastic partial differential equation;
colored noise; boundary conditions;
non-Lipschitz coefficients; pathwise uniqueness;
comparison theorem; H\"older continuity.
}}
\section{Introduction and main results}
\setcounter{equation}{0}

It is first shown in Konno and Shiga \cite{KoS88} (also in Reimers \cite{Rei89})
that the one-dimensional super-Brownian motion with binary
branching has a
jointly continuous density
$\{X_t(x):t>0,x\in\mbb{R}\}$ which satisfies the
following stochastic partial differential equation
(SPDE):
 \beqlb\label{1.0}
\frac{\partial}{\partial t} X_t(x)
 =
\frac{1}{2}\Delta X_t(x) + X_t(x)^\beta\dot{W}_t(x), \quad
t>0,
~x\in\mbb{R},
 \eeqlb
where $\Delta$ denotes the one-dimensional Laplacian operator and
$\dot{W}$ is a
space-time Gaussian white noise.
The weak uniqueness of the solution to \eqref{1.0}
follows from the uniqueness of solution to the martingale problem for
 the associated
super-Brownian motion.
We refer to \cite{D,dynkin,Eth,Li11,Per} for introduction to super-Brownian motion.
The pathwise uniqueness of the nonnegative solution to \eqref{1.0}
remains open even though it is studied by many authors;
see \cite{Dawsonli06,Dawsonli12,MyP11,Xio13} and the references
therein.
For an one-dimensional super-Brownian motion in random environment,
its density process satisfies more general SPDE than equation \eqref{1.0}
(see Dawson \textit{et al.} \cite{Dawson}),
and the joint H\"older continuity is studied
in Li \textit{et al.} \cite{LiWXZ} and the improved result is obtained
in Hu \textit{et al.} \cite{HuLD}.
If the noise $\dot{W}$ is colored
in space and white in time,
the existence and pathwise uniqueness of nonnegative
solution to the SPDE were established
by Sturm \cite{Sturm03} and Mytnik \textit{et al.}
\cite{MPS06}, respectively.
Further work have been studied in Rippl and Sturm \cite{RS13} and
Neuman \cite{N14}.
For certain diffusion coefficient depending on the
spatial derivative of the solution,
Gomez \textit{et al.} \cite{Gomez13} and Xiong and Yang \cite{XiongYang17a} studied the pathwise uniqueness and strong existence
of the solution, respectively. We also refer the reader to Xiong \cite{X} for some other related SPDEs connected to superprocesses.

There are many authors studying the following SPDE where the spatial domain is a finite interval, $[0,1]$, with
Dirichlet or Neumann conditions imposed on the endpoints:
 \beqlb\label{1.1}
\frac{\partial Y_t(x)}{\partial t}=\frac12\Delta Y_t(x)
+G(Y_t(x))
+H(Y_t(x))\dot{W}_t(x),
 \eeqlb
where $x\in(0,1)$, $t\ge0$, $G$ and $H$ are continuous functions on $\mbb{R}$,
and $\dot{W}$ denote the
space-time Gaussian white noise on $[0,\infty)\times[0,1]$.
It is shown in \cite[Chapter 3]{Walsh} that SPDE \eqref{1.1} has a unique solution under Lipschitz continuity
on $H$ and $G$.
If $H$ is a constant and $G$ is
locally bounded and satisfies one sided linear growth condition,
then SPDE \eqref{1.1} has a unique strong solution and the comparison theorem
holds, which is given in \cite{GyPardoux93a} and
extended to nonnecessarily locally bounded for $G$
in \cite{GyPardoux93b}.
If the assumptions concerning $G$ are those of \cite{GyPardoux93a},
and $H$ has locally Lipschitz derivative and satisfies
a linear growth condition, then the existence and uniqueness and
comparison theorem
are studied in \cite{BallyGyPardoux94} and the results are extended
to Lipschitz continuous $G$ in \cite{Gy95}.
There are also many authors concerning on
white noise driven SPDEs with reflection,
which adds a certain random measure on
the right hand of \eqref{1.1} and
the existence and uniqueness of solution and comparison theorem
are studied in \cite{Donati-MartinPardoux93,NualartPardoux92,XuZhang09,Zhang11} with $G$ and $H$ satisfying Lipschitz and linear growth conditions.
Various properties of the solution were studied in
\cite{DalangMueller06,DalangZhang13,Kalsi19,YangZhang14,Zhang10}.

In this paper we study SPDE \eqref{1.1} with
Dirichlet, Neumann or mixed {\em nonhomogeneous random boundary conditions},
where $\dot{W}$ is a Gaussian  noise that is white in time and colored in
space. Namely,
it is a Gaussian martingale measure on $[0,\infty)\times[0,1]$
in the sense of \cite[Chapter 2]{Walsh} and can be characterized by its covariance functional
 \beqnn
\mbf{E}\big[W(\phi)W(\psi)\big]
:=\int_0^\infty\int_0^1\int_0^1 \phi(s,x)\psi(s,y)\kappa(x,y)\dd s\dd x\dd y
 \eeqnn
for each continuous functions $\phi,\psi$ on $[0,\infty)\times[0,1]$,
where $\kappa$ is a nonnegative bounded function on $[0,1]\times[0,1]$ and $\kappa(x,y)=\kappa(y,x)$ for all $x,y\in[0,1]$.

To continue with the introduction we state some notation.
Let $\mathscr{B}[0,1]$ denote the set of Borel functions on $[0,1]$.
For $f,g\in \mathscr{B}[0,1]$ let $\<f,g\>=\int_0^1f(x)g(x)\dd x$ whenever it exists.
Let $B[0,1]$ be the Banach space of bounded measurable functions
on $[0,1]$ furnished with the supremum norm $\|\cdot\|_0$.
We use $C[0,1]$ to denote the subset of continuous functions on $[0,1]$.
Define $C(\mbb{R})$, $C([0,\infty)\times[0,1])$ and $C([0,T]\times[0,1])$ similarly.
For any integer $n\ge1$ let $C^n[0,1]$ be the subset of
$C[0,1]$ of functions with bounded continuous derivatives up to the $n$th order.
Let $C_c^n[0,1]$ be the subset of $C^n[0,1]$ of functions with the supports in $(0,1)$.

Suppose that
$(\mu_0(t))_{t\ge0}$ and $(\mu_1(t))_{t\ge0}$ are two continuous processes
and that $Y_0\in C[0,1]$.
In this paper we aim to prove the strong existence, pathwise uniqueness
and to establish a comparison theorem for the solution to SPDE \eqref{1.1} with $G$ satisfying certain non-Lipschitz
condition and $H$ satisfying H\"older condition and
with one of the following boundary conditions:
\begin{itemize}
\item[{\normalfont(C1)}]
$Y_t(0)=\mu_0(t)$
and $Y_t(1)=\mu_1(t)$ almost surely for all $t\ge0$;
\item[{\normalfont(C2)}]
$\nabla Y_t(0)=\mu_0(t)$
and $\nabla Y_t(1)=\mu_1(t)$ almost surely for all $t\ge0$;
\item[{\normalfont(C3)}]
$Y_t(0)=\mu_0(t)$
and $\nabla Y_t(1)=\mu_1(t)$ almost surely for all $t\ge0$;
\item[{\normalfont(C4)}]
$\nabla Y_t(0)=\mu_0(t)$
and $Y_t(1)=\mu_1(t)$ almost surely for all $t\ge0$.
\end{itemize}

To present the definition of the solution to SPDE (\ref{1.1}) with various boundary conditions precisely, we introduce the corresponding boundary conditions on the test functions
 $f\in C^2[0,1]$:
\begin{itemize}
\item[{\normalfont(D1)}]
$f(0)=f(1)=0$;
\item[{\normalfont(D2)}]
$f'(0)=f'(1)=0$;
\item[{\normalfont(D3)}]
$f(0)=f'(1)=0$;
\item[{\normalfont(D4)}]
$f'(0)=f(1)=0$.
\end{itemize}

\bdefinition\label{d1.1}
Let $i\in\{1,2,3,4\}$.
We say that $(Y_t)_{t\ge0}$ on a filtered probability space is a weak solution to \eqref{1.1} with boundary condition (Ci) if
the mapping $(t,x)\mapsto Y_t(x)$ is in $C([0,\infty)\times[0,1])$ almost surely and
there exist a Gaussian colored noise $W$, an initial $Y_0$ and
boundary conditions $\mu_0(t)$, $\mu_1(t)$, so that
for any $f\in C^2[0,1]$
satisfying (Di), we have
 \beqlb\label{1.001}
\<Y_t,f\>
 \ar=\ar
\<Y_0,f\>+\frac12\int_0^t[\<Y_s,f''\>+F^{(i)}_s(f)]\dd s \cr
 \ar\ar
+\int_0^t\<G(Y_s),f\>\dd s
+\int_0^t\int_0^1H(Y_s(x))f(x)W(\dd s,\dd x),\quad
t\ge0
 \eeqlb
almost surely, where
 \beqlb\label{1.6}
F^{(1)}_s(f):=f'(0)\mu_0(s)-f'(1)\mu_1(s),\quad
F^{(2)}_s(f):=-f(0)\mu_0(s)+f(1)\mu_1(s)
 \eeqlb
and
 \beqlb\label{1.7}
F^{(3)}_s(f):=f'(0)\mu_0(s)+f(1)\mu_1(s),\quad
F^{(4)}_s(f):=-f(0)\mu_0(s)-f'(1)\mu_1(s).
 \eeqlb

SPDE \eqref{1.1} with boundary condition (Ci) has a strong solution if for any Gaussian colored noise $W$ on filtered probability space
$(\Omega,\mathscr{F},\mathscr{F}_t,\mbf{P})$, for given initial $Y_0$ and
the continuous processes $(\mu_0(t))_{t\ge0}$ and
$(\mu_1(t))_{t\ge0}$,
there exists a process $(Y_t)_{t\ge0}$
 so that
$(t,x)\mapsto Y_t(x)$ is in $C([0,\infty)\times[0,1])$
and \eqref{1.001} holds for all $f\in C^2[0,1]$
satisfying (Di).
\edefinition

We also study the joint H\"older continuity of the solution. We aim to prove that
the H\"older exponent in spatial variable is arbitrarily
close to 1/2 and in time variable is arbitrarily
close to 1/4 for the case of boundary condition (C2); the H\"older exponents of the solution in spatial
and in time variables depend on the H\"older exponent of the non-derivative boundary condition
for the case of
boundary conditions (C1), (C3) and (C4).
Throughout this paper we always assume that all
random variables defined on the same filtered probability space
$(\Omega,\mathscr{F},\mathscr{F}_t,\mbf{P})$.
Let $\mbf{E}$ denote the corresponding expectation.

The following conditions will be imposed in the rest of the paper.

\bcondition\label{tc}
There are constants $\gamma\in[2^{-1},1]$ and $C>0$ so that
 \beqlb\label{1.2}
|G(x)|\le C[|x|+1],\qquad x\in\mbb{R}
 \eeqlb
and
 \beqlb\label{1.3}
|H(x)-H(y)|\le C|x-y|^\gamma,\qquad x,y\in\mbb{R}.
 \eeqlb
\econdition

\bcondition\label{tc2}
Let $r_0>0$ be a constant.
For each $n\ge1$ there is a non-decreasing and concave function
$r_n$ on $[0,\infty)$ so that $\int_{0+}\frac{1}{r_n(x)}\dd x=\infty$
and
 \beqnn
(x-y)r_0+G(x)-G(y)\le r_n(x-y),\qquad -n\le y\le x\le n.
 \eeqnn
\econdition

\btheorem\label{t1.1} (Comparison theorem)
Suppose that $G_1$ and $G_2$ are two continuous functions
satisfying $G_1(y)\le G_2(y)$ for all $y\in\mbb{R}$, and that Condition \ref{tc2} holds for $G_2$
and \eqref{1.3} holds for $H$.
For each $i=1,2$ let $(Y_t^{(i)})_{t\ge0}$ be a weak solution to \eqref{1.1}
with $G$ replaced by $G_i$.
We also assume that $(Y_t^{(1)})_{t\ge0}$ and $(Y_t^{(2)})_{t\ge0}$
satisfy the same boundary condition.
If $Y_0^{(1)}(x)\le Y_0^{(2)}(x)$ for all $x\in[0,1]$,
then $\mbf{P}\{Y_t^{(1)}(x)\le Y_t^{(2)}(x)$ for all $t\ge0$ and $x\in[0,1]\}=1$.
\etheorem
Applying this theorem we obtain the pathwise uniqueness to \eqref{1.1} immediately.
\bcorollary\label{t1.3} (Pathwise uniqueness)
Suppose that \eqref{1.3} and Condition \ref{tc2} hold.
Let $(Y_t^{(1)})_{t\ge0}$ and $(Y_t^{(2)})_{t\ge0}$ be two weak solutions to \eqref{1.1} with the same boundary condition.
If $Y_0^{(1)}(x)= Y_0^{(2)}(x)$ for all $x\in[0,1]$,
then $\mbf{P}\{Y_t^{(1)}(x)= Y_t^{(2)}(x)$ for all $t\ge0$ and $x\in[0,1]\}=1$.
\ecorollary

\btheorem\label{t1.2} (Existence)
Suppose that Conditions \ref{tc} and \ref{tc2} hold.
Then \eqref{1.1} has a unique strong solution $(Y_t)_{t\ge0}$
with one of the boundary conditions (C1) to (C4).
Moreover, if $Y_0\ge0$, $G(0)=H(0)=0$ and $\mu_0(t)=\mu_1(t)=0$ for all $t\ge0$, then $\mbf{P}\{Y_t\ge0$ for all $t\ge0\}=1$.
\etheorem

The next result is concerned with
 the H\"older continuity of the
solution both in time  and in space variables.
We first state the following assumption.

\bcondition\label{tc4.1}
(i) There is a constant $C>0$ so that
 \beqnn
|G(x)|\le C(|x|+1),\quad |H(x)|\le C(|x|+1),\qquad x\in\mbb{R}.
 \eeqnn

(ii) For each $T>0$ and $p>1$,
 \beqnn
\mbf{E}\Big[\sup_{0<t\le T}[|\mu_0(t)|^{2p}+|\mu_1(t)|^{2p}]\Big]<\infty.
 \eeqnn

(iii) There are constants $0<\gamma_0<1/2$
and $C>0$ so that
\begin{itemize}
\item
For boundary condition (C1),
 \beqlb\label{7.13}
\mbf{E}\Big[|\mu_i(t)-\mu_i(s)|^p\Big]\le C|t-s|^{p\gamma_0},~~
0<t,s\le T
 \eeqlb
holds for $i=1,2$;
\item
For boundary condition (C3), \eqref{7.13} holds for $i=0$;
\item
For boundary condition (C4), \eqref{7.13} holds for $i=1$.
\end{itemize}
\econdition

\btheorem\label{t6.1}
(Joint H\"older continuity)
Suppose that $(Y_t)_{t\ge0}$ is a weak solution
to \eqref{1.1}
with one of the boundary conditions (C1) to (C4) and that Condition \ref{tc4.1} holds
and $0<T_1<T$.
Then $[T_1,T]\times[0,1]\ni(t,x)\mapsto Y_t(x)$ is H\"older
continuous with exponent $\eta_1/2$ in the time variable and with exponent $\eta_2$ in
space variables,
where $\eta_1,\eta_2\in(0,1/2)$ in the case of boundary condition (C2)
and $\eta_1,\eta_2\in(0,\gamma_0)$ in the case of boundary conditions (C1), (C3) and (C4).
Namely, there exists a random variable $K\ge0$ depending on $\eta_1$ and $\eta_1$ so that
 \beqnn
|Y_{t_1}(x_1)-Y_{t_2}(x_2)|
\le K(|t_1-t_2|^{\eta_1/2}+|x_1-x_2|^{\eta_2}),\qquad
t_1,t_2\in[T_1,T],~~x_1,x_2\in[0,1].
 \eeqnn
\etheorem

In the following we consider SPDE \eqref{1.1}
on real semi-axis.

\btheorem\label{t1.2'}
Suppose that Conditions \ref{tc} and \ref{tc2} hold.
Then there is a unique strong solution $(Y_t)_{t\ge0}$
to \eqref{1.1}
with $x\in(-\infty,0)$ (or $x\in(0,\infty)$) and one of the following
boundary conditions:
\begin{itemize}
\item[{\normalfont(C1')}]
$Y_t(0)=\mu_0(t)$ almost surely for all $t\ge0$;
\item[{\normalfont(C2')}]
$\nabla Y_t(0)=\mu_0(t)$ almost surely for all $t\ge0$.
\end{itemize}
Moreover, if $Y_0\ge0$, $G(0)=H(0)=0$ and $\mu_0(t)=\mu_1(t)=0$ for all $t\ge0$, then $\mbf{P}\{Y_t\ge0$ for all $t\ge0\}=1$.
\etheorem

The rest of the paper is organized as follows.
In Section 2, we state some properties of the solution
to heat equations with boundary conditions,
which are used in the proofs of theorems.
Theorem \ref{t1.1} is proved in Section 3.
In Section 4, we establish the proof of Theorems \ref{t1.2} and \ref{t1.2'}.
We derive the proofs of Theorem \ref{t6.1} in Section 5.

\textit{Notation:}
Let $\nabla$ and $\Delta$ be the first and the second order spatial differential operators, respectively.
Similarly, $\nabla_x$ and $\Delta_x$ denote the first and the second order spatial
differential operators with respect to the variable $x$.
Let $\partial_x$ denote the first partial derivative with respect to the variable $x$ and $\kappa_0:=\sup_{y_1,y_2\in[0,1]}\kappa(y_1,y_2)$.
We use $C$ to denote a positive constant whose value might change from line to line.

\section{Heat equations with boundary conditions}

\setcounter{equation}{0}

In this section we establish some properties of the solutions
to the following heat equations with different boundary conditions:
\begin{equation}\label{1.4}
  \left\{
   \begin{aligned}
   \partial_tp_t(x,y) &= \frac12\Delta_x p_t(x,y),\quad t>0,~x,y\in(0,1),\\
   \lim_{t\downarrow0}p_t(x,y) &=\delta_y(x),\qquad\qquad x,y\in[0,1].
   \end{aligned}
   \right.
  \end{equation}
For $t>0$ and $x,y\in\mbb{R}$ set $q_t(x):=\frac{1}{\sqrt{2\pi t}}\exp\{-x^2/(2t)\}$
and $q_t(x,y):=q_t(x-y)$.
For each $t>0$ and $x,y\in\mbb{R}$ let
 \begin{gather}
\label{2.1}p_t^{(1)}(x,y)
 =
\sum_{k=-\infty}^\infty[q_t(2k+x,y)-q_t(2k-x,y)], \\
\label{2.2}
p_t^{(2)}(x,y)
 =
\sum_{k=-\infty}^\infty[q_t(2k+x,y)+q_t(2k-x,y)]
 \end{gather}
and
 \begin{gather}
\label{2.4}p_t^{(3)}(x,y)
 =
2\sum_{k=-\infty}^\infty[q_t(4k+x,y)-q_t(4k-x,y)]
-\sum_{k=-\infty}^\infty[q_t(2k+x,y)-q_t(2k-x,y)], \\
\label{2.3}p_t^{(4)}(x,y)
 =
2\sum_{k=-\infty}^\infty[q_t(4k+x,y)+q_t(4k-x,y)]
-\sum_{k=-\infty}^\infty[q_t(2k+x,y)+q_t(2k-x,y)].
 \end{gather}
It is elementary to check that
for each $i=1,2,3,4$, $p_t^{(i)}(x,y)$ is the unique solution to \eqref{1.4} with $p_t^{(i)}(x,\cdot)$ satisfying the
boundary condition (Di) for all $t>0$ and $x\in\mbb{R}$.

For simplicity we write $p_t(x,y)$ for $p_t^{(i)}(x,y)$ when the
conclusion is true for all $i$. We write $P_tf(y)=\int_0^1p_t(x,y)f(x)\dd x$ for $x\in[0,1]$
and $f\in\mathscr{B}[0,1]$.
It is elementary to check the following two lemmas.
\blemma\label{t2.5}
For each $t\ge0$ and $f\in C[0,1]$
we have
 \beqnn
p_t(x,y)=p_t(y,x),\quad
\int_0^1|p_t(x,y)|\dd y
+\int_0^1|p_t(x,y)|\dd x\le 12,\qquad x,y\in\mbb{R}.
 \eeqnn
and
 \beqnn
\lim_{t\to0}\int_0^1f(x)p_t(x,y)\dd x=f(y),\quad
\lim_{t\to0}\int_0^1f(y)p_t(x,y)\dd y=f(x),\qquad x,y\in\mbb{R}.
 \eeqnn
\elemma

\blemma\label{t2.1}
For each $t,s\ge0$ and $f\in C[0,1]$ we have
 \beqnn
P_tP_sf(x)=P_{t+s}f(x),\qquad t,s>0.
 \eeqnn
\elemma

\blemma\label{t2.2}
Let $T>0$ be fixed.
Then
 \beqnn
|p_t(x_1,y)-p_t(x_2,y)|\le
Ct^{-1}|x_1-x_2|,\qquad t\in(0,T],~x_1,x_2,y\in[0,1]
 \eeqnn
and
 \beqnn
|p_t(x,y)|\le Ct^{-1/2},~
|p_{t+\varepsilon}(x,y)-p_t(x,y)|\le
Ct^{-1}\varepsilon^{1/2},\quad t\in(0,T],~\varepsilon>0,~x,y\in[0,1].
 \eeqnn
\elemma
\proof
It is obvious that
 \beqlb\label{2.5}
\sum_{k=-\infty}^{\infty}q_t(k+x-y)\le C t^{-1/2},\qquad x,y\in[0,1],~t\in(0,T].
 \eeqlb
It then follows from \eqref{2.1}-\eqref{2.3} that
 \beqnn
|p_t(x,y)|\le Ct^{-1/2},\qquad x,y\in[0,1],~t\in(0,T].
 \eeqnn
Due to \cite[(2.4e)]{Rosen87} we obtain
 \beqlb\label{4.8}
|q_t(x_1,y)-q_t(x_2,y)|
\le
Ct^{-1/2}|x_1-x_2|[q_{4t}(x_1-y)+q_{4t}(x_2-y)]
 \eeqlb
for $x_1,x_2,y\in\mbb{R}$ and $t>0$.
From \eqref{2.1}-\eqref{2.3} it follows that for $x_1,x_2,y\in[0,1]$ and $t\in(0,T]$,
 \beqnn
 \ar\ar
|p_t(x_1,y)-p_t(x_2,y)| \cr
 \ar\le\ar
Ct^{-1/2}|x_1-x_2|\sum_{k=-\infty}^{+\infty}
\big[q_{4t}(4k+x_1,y)+q_{4t}(4k-x_2,y) \cr
 \ar\ar\quad\qquad\qquad\qquad\qquad\quad
+q_{4t}(2k+x_1,y)+q_{4t}(2k-x_2,y)\big] \cr
 \ar\le\ar
Ct^{-1}|x_1-x_2|,
 \eeqnn
which gives the first assertion \eqref{2.5}.
Observe that $q_t(x)=t^{-1/2}q_1(t^{-1/2}x)$
and  $xq_t(x)\le C$ for all $t,x>0$.
From \eqref{4.8} it follows that for $0<t\le T$ and $x>0$,
 \beqnn
 \ar\ar
|q_{t+\varepsilon}(x)
-q_{t}(x)| \cr
 \ar\le\ar
|(t+\varepsilon)^{-1/2}-t^{-1/2}|q_1((t+\varepsilon)^{-1/2}x)
+
t^{-1/2}|q_1((t+\varepsilon)^{-1/2}x)-q_1(t^{-1/2}x)| \cr
 \ar\le\ar
\varepsilon^{1/2}t^{-1/2}(t+\varepsilon)^{-1/2}q_1((t+\varepsilon)^{-1/2}x) \cr
 \ar\ar
+
Ct^{-1/2}x|(t+\varepsilon)^{-1/2}-t^{-1/2}|\cdot [q_{4}((t+\varepsilon)^{-1/2}x)+q_{4}(t^{-1/2}x)] \cr
 \ar\le\ar
\varepsilon^{1/2}t^{-1}q_1((T+\varepsilon)^{-1/2}x)
+\varepsilon^{1/2}t^{-1}C[q_8((T+\varepsilon)^{-1/2}x)+q_8(T^{-1/2}x)]
\cr
 \ar\le\ar
C\varepsilon^{1/2}t^{-1}[(T+\varepsilon)^{-1/2}q_{8(T+\varepsilon)}(x)
+T^{-1/2}q_{8T}(x)].
 \eeqnn
Combining this with \eqref{2.1}-\eqref{2.3} we obtain the second assertion.
\qed

\blemma\label{t5.2}
We have
 \beqnn
\sup_{0<t\le T,\,x\in[0,1]}\sum_{n=-\infty}^\infty \int_0^t|\nabla q_s(n-x)|\dd s<\infty.
 \eeqnn
\elemma
\proof
Note that for $v>0$,
 \beqlb\label{7.6}
|\nabla q_v(1)|\le C(1_{\{v\le1\}}+v^{-3/2}1_{\{v>1\}}).
 \eeqlb
It then follows from a change of variable that
for each $y>0$,
 \beqlb\label{7.7}
\int_0^t|\nabla q_s(y)|\dd s
=\int_0^{ty^{-2}}|\nabla q_v(1)|\dd v
\le\int_0^\infty |\nabla q_v(1)|\dd v<\infty.
 \eeqlb

For all $y>1$,
 \beqnn
|\nabla q_s(y)|
\le
y^2s^{-1}q_s(y)
\le Cq_{2s}(y),
 \eeqnn
which leads to
 \beqnn
\sum_{|n|\ge 2}\int_0^t|\nabla q_s(n-x)|\dd s
\le
C\int_0^t\sum_{|n|\ge 2}|q_{2s}(n-x)|\dd s
\le
Ct,\quad x\in[0,1].
 \eeqnn
Combining the above inequality with \eqref{7.7} one ends the proof.
\qed

\blemma\label{t5.3}
Let $p\ge1$ be fixed and $(h(t))_{t\ge0}$ be a stochastic process
satisfying
 \beqnn
\mbf{E}\Big[\sup_{0<t\le T}|h(t)|^p\Big]<\infty.
 \eeqnn
Suppose that there is a constant $0<\tilde{\gamma}_0<1/2$ so that
 \beqnn
\mbf{E}\big[|h(s+t)-h(t)|^p\big]\le Cs^{p\tilde{\gamma}_0},\qquad s,t>0.
 \eeqnn
Then
 \beqnn
\mbf{E}\Big[\Big|\int_0^th(s)\sum_{n=-\infty}^\infty
[\nabla q_s(2n-x_1)
-\nabla q_s(2n-x_2)]\dd s\Big|^p
\Big]
\le
C[t^{-p/2}+1]|x_1-x_2|^{p\tilde{\gamma}_0}
 \eeqnn
for all $t\in(0,T]$ and $x_1,x_2\in[0,1]$.
\elemma
\proof
Let
\[I_t(x,y):=\int_0^th(s)\nabla q_s(x+y)\dd s-\int_0^th(s)\nabla q_s(y)\dd s,\qquad 0\le x,\ y\le 1.\]
We first show that for $0<x,y\le 1$,
 \beqlb\label{7.8}
\mbf{E}\left[\left|I_t(x,y)\right|^p
\right]
\le
C[t^{-p/2}+1]x^{p\tilde{\gamma}_0}.
 \eeqlb
By a change of variable and \eqref{7.6},
 \beqlb\label{7.4}
|I_t(x,y)|
 \ar=\ar
\Big|\int_0^{t(x+y)^{-2}}h(v(x+y)^2)\nabla q_v(1)\dd v
-\int_0^{ty^{-2}}h(vy^2)\nabla q_v(1)\dd v\Big| \cr
 \ar\le\ar
\int^{ty^{-2}}_{t(x+y)^{-2}}|h(vy^2)\nabla q_v(1)|\dd v
+I_{t,1}(x,y) \cr
 \ar\le\ar
C\sup_{0<v\le T}|h(v)|t^{-1/2}x
+I_{t,1}(x,y),
 \eeqlb
where
\[I_{t,1}(x,y)=\int_0^{t(x+y)^{-2}}|[h(v(x+y)^2)-h(vy^2)]\nabla q_v(1)|\dd v.\]
Taking $\delta>0$ satisfying $\tilde{\gamma}_0+1/p<3\delta/2<1/2+1/p$, then
 \beqnn
\frac{3(1-\delta)p}{2(p-1)}>1,\quad
\frac{3p\delta}{2}-p\tilde{\gamma}_0>1.
 \eeqnn
By H\"older's inequality and \eqref{7.6} we find
 \beqnn
 \ar\ar
\mbf{E}\big[|I_{t,1}(x,y)|^p\big] \cr
 \ar\le\ar
\Big|\int_0^\infty|\nabla q_v(1)|^{\frac{(1-\delta)p}{(p-1)}}\dd v\Big|^{p-1}
\cdot\int_0^\infty\mbf{E}\big[|h(v(x+y)^2)-h(vy^2)|^p\big]
\cdot|\nabla q_v(1)|^{p\delta}\dd v \cr
 \ar\le\ar
Cx^{p\tilde{\gamma}_0}
\Big|\int_0^{\infty}
\big[1_{\{v\le 1\}}+v^{-\frac{3(1-\delta)p}{2(p-1)}}1_{\{v> 1\}}\big]\dd v\Big|^{p-1}\int_0^\infty
\big[1_{\{v\le 1\}}+v^{p\tilde{\gamma}_0-\frac{3p\delta}{2}}1_{\{v>1\}}\big]\dd v\cr
\ar\le\ar
C x^{p\tilde{\gamma}_0}
 \eeqnn
for all $0<x,y\le1$ and $t>0$,
which gives \eqref{7.8} by \eqref{7.4}.

Observe that
 \beqnn
|\Delta q_t(x)|\le Cq_{2t}(x),\qquad t>0,~x\ge1.
 \eeqnn
Then, for all $x,t>0$ and $y\ge1$, we have
 \beqnn
|\nabla q_t(x+y)-\nabla q_t(y)|
\le
x\int_0^1|\Delta q_t(y+x\theta)|\dd \theta
\le
Cx\int_0^1q_{2t}(y+x\theta)\dd \theta
\le
Cxq_{2t}(y),
 \eeqnn
which implies that for all $x,y,t>0$ and $x+y\le1$,
 \beqnn
 \ar\ar
\sum_{|n|\ge1}\int_0^t|\nabla q_s(2n-x-y)-\nabla q_s(2n-y)|\dd s \cr
 \ar\ar\quad
\le Cx\int_0^t\Big[\sum_{n\le-1}q_{2s}(2n-y)+\sum_{n\ge1}q_{2s}(2n-1-y)\Big]\dd s
\le Cxt.
 \eeqnn
Combining this with \eqref{7.8} one ends the proof.
\qed

\section{Proof of Theorem \ref{t1.1}}
\setcounter{equation}{0}
In this section we prove Theorem \ref{t1.1}
by a modification of the Yamada-Watanabe argument for ordinary stochastic
differential equations.

\noindent{\it Proof of Theorem \ref{t1.1}}.
For $k\ge1$ put $a_k=\exp\{-k(k+1)/2\}$. Let $\psi_k\in
C_c^2[0,1]$ satisfy $\mbox{supp}(\psi_k)\subset(a_k,a_{k-1})$,
$\int_{a_k}^{a_{k-1}}\psi_k(x)dx=1$, and $0\le\psi_k(x)\le 2/(kx)$
for all $x>0$ and $k\ge1$.
Let $x^+:=x\vee0$ and
$\phi_k(x)=\int_0^{x^+}dy\int_0^y\psi_k(z)\dd z$ for $x\in\mbb{R}$ and
$k\ge1$. Then for all $x\in\mbb{R}$, $0\le\phi'_k(x)\le 1$,
$\phi_k'(x)\to 1_{\{x>0\}}$ and $\phi_k(x)\to x^+$ as
$k\to\infty$.

For $x\in{\mbb{R}}$ define $J(x)=\int_{\mbb{R}}\e^{-|y|}\rho_0(x-y)dy$
with the mollifier $\rho_0$ given by
 \beqnn
\rho_0(x)=c_0\exp\big(-1/(1-x^2)\big)1_{\{|x|<1\}},
 \eeqnn
where $c_0>0$ is a constant so that $\int_{\mbb{R}}\rho_0(x)dx=1$.
Moreover, due to (2.1) of \cite{Mitoma}, for each $n\ge0$ there exist constants
$\bar{C}_n,\tilde{C}_n>0$ so that
 \beqlb\label{5.16a}
\bar{C}_n\e^{-|x|}\le |J^{(n)}(x)|\le \tilde{C}_n\e^{-|x|}, \qquad
x\in\mbb{R}.
 \eeqlb
For $\zeta>0$ and $x\in\mbb{R}$ let $J_\zeta(x)=(J(x))^\zeta$.
Then by using \eqref{5.16a} one can see that there is a constant $C_0>0$ independent of $\zeta$ so that
 \beqlb\label{5.16}
|J_\zeta''(x)|\le C_0\zeta(\zeta+1)J_\zeta(x),
\qquad \zeta>0,~x\in\mbb{R}.
 \eeqlb

For the solutions $(Y_t^{(1)})_{t\ge0}$ and $(Y_t^{(2)})_{t\ge0}$ to \eqref{1.1}
with boundary condition (Ci),
let $p_t(x,y)$ denote the solution to
\eqref{1.4} with the boundary condition (Di).
It then follows that for $k=1,2$,
 \beqlb\label{3.5}
\<Y_t^{(k)},p_\delta^x\>
 \ar=\ar
\<Y_0^{(k)},p_\delta^x\>+\frac12\int_0^t\big[\Delta_x\<Y_s^{(k)},p_\delta^x\>
+F_s^{(i)}(p_\delta^x)]\dd s \cr
 \ar\ar
+\int_0^t\<G(Y_s^{(k)}),p_\delta^x\>\dd s
+\int_0^t\int_0^1 H(Y_s^{(k)}(y))p_\delta^x(y) W(\dd s,\dd y),
 \eeqlb
where $p^x_\de=p_\de(x,\cdot)$, $F_s^{(i)}(p_\delta^x)$ is defined in \eqref{1.6} and \eqref{1.7}.
For each $t\ge0$ and $k=1,2$ we extend the definition of $Y_t^{(k)}$
on $\mbb{R}$ by $Y_t^{(k)}(x)=Y_t^{(k)}(0)$ for $x\le0$ and
$Y_t^{(k)}(x)=Y_t^{(k)}(1)$ for $x\ge1$.

For $n\ge1$, we define stopping time
 \beqnn
\tau_n:=\inf\big\{t\ge0:\|Y_t^{(1)}\|_0+\|Y_t^{(2)}\|_0\ge n\big\}.
 \eeqnn
Set $v_t(x)=Y_t^{(1)}(x)-Y_t^{(2)}(x)$ and $v_t^\delta(x)=\<v_t,p_\delta^x\>$.
From \eqref{3.5} it follows that
 \beqlb\label{5.17}
v_{t\wedge\tau_n}^\delta(x)
 \ar=\ar
v_0^\delta(x)+
\frac12\int_0^{t\wedge\tau_n}\Delta_x v_s^\delta(x)\dd s
+\int_0^{t\wedge\tau_n}[G^{\delta,1}(x)+G^{\delta,2}(x)]\dd s \cr
 \ar\ar
+\int_0^{t\wedge\tau_n}\int_0^1M_s^\delta(x,y)W(\dd s,\dd y),
 \eeqlb
where
 \beqnn
G^{\delta,1}(x):=\<G_1(Y_s^{(1)})-G_2(Y_s^{(1)}),p_\delta^x\>,~~
G^{\delta,2}(x):=\<G_2(Y_s^{(1)})-G_2(Y_s^{(2)}),p_\delta^x\>
 \eeqnn
and
 \beqnn
M_s^\delta(x,y):=[H(Y_s^{(1)}(y))-H(Y_s^{(2)}(y))]p_\delta^x(y).
 \eeqnn
By \eqref{1.3} and H\"older's inequality we obtain
 \beqlb\label{5.6}
M^\delta_s(x)
 \ar:=\ar
\Big|\int_0^1\dd y_1\int_0^1
M_s^\delta(x,y_1)
M_s^\delta(x,y_2)\kappa(y_1,y_2)\dd y_2\Big| \cr
 \ar\le\ar
\kappa_0\Big[\int_0^1
\big|[H(Y_s^{(1)}(y))-H(Y_s^{(2)}(y))]p_\delta^x(y)\big| \dd y\Big]^2 \cr
 \ar\le\ar
C \Big[\int_0^1
|v_s(y)|^{\gamma}\cdot|p_\delta^x(y)| \dd y\Big]^2
\le
C\<|v_s|,|p_\delta^x|\>^{2\gamma}.
 \eeqlb
It then follows from \eqref{5.17} and It\^o's formula that
 \beqnn
\phi_k(v_{t\wedge\tau_n}^\delta(x))
 \ar=\ar
\phi_k(v_0^\delta(x))
+
\frac12\int_0^{t\wedge\tau_n}\big[\phi'_k(v_s^\delta(x))
\Delta_x v_s^\delta(x)+\phi''_k(v_s^\delta(x))M^\delta_s(x)\big]\dd s \cr
 \ar\ar
+\int_0^{t\wedge\tau_n}\phi'_k(v_s^\delta(x))[G^{\delta,1}(x)+G^{\delta,2}(x)]\dd s \cr
 \ar\ar
+\int_0^{t\wedge\tau_n}\int_0^1\phi'_k(v_s^\delta(x))M_s^\delta(x,y)W(\dd s,\dd y),
 \eeqnn
which leads to
 \beqlb\label{5.5}
  \ar\ar
\int_{\mbb{R}}\mbf{E}\big[\phi_k(v_{t\wedge\tau_n}^\delta(x))\big]J_\zeta(x)\dd x
-\int_{\mbb{R}}\phi_k(v_0^\delta(x))J_\zeta(x)\dd x\cr
 \ar=\ar
\mbf{E}\Big[\int_0^{t\wedge\tau_n}\Big[2^{-1}\int_{\mbb{R}}\phi'_k(v_s^\delta(x))
\Delta_x v_s^\delta(x)J_\zeta(x)\dd x
+2^{-1}\int_{\mbb{R}}\phi''_k(v_s^\delta(x))M^\delta_s(x)J_\zeta(x)\dd x \cr
 \ar\ar
+\int_{\mbb{R}}\phi'_k(v_s^\delta(x))G^{\delta,1}(x)J_\zeta(x)\dd x
+\int_{\mbb{R}}\phi'_k(v_s^\delta(x))G^{\delta,2}(x)J_\zeta(x)\dd x\Big]\dd s \cr
 \ar=:\ar
\mbf{E}\Big[\int_0^{t\wedge\tau_n}\big[2^{-1}I_{1,k}^\delta(s)
+2^{-1}I_{2,k}^\delta(s)
+I_{3,k}^\delta(s)+I_{4,k}^\delta(s)\big]\dd s\Big].
 \eeqlb
By integration by parts and \eqref{5.16},
 \beqnn
I_{1,k}^\delta(s)
 \ar=\ar
\int_{\mbb{R}}\Delta_x(\phi_k(v_s^\delta(x)))J_\zeta(x)\dd x
-\int_{\mbb{R}}\phi''_k(v_s^\delta(x))|\nabla_x
v_s^\delta(x)|^2J_\zeta(x)\dd x \cr
 \ar\le\ar
\int_{\mbb{R}}\phi_k(v_s^\delta(x))J_\zeta''(x)\dd x
\le
C_1\zeta(\zeta+1)\int_{\mbb{R}}
\phi_k(v_s^\delta(x))J_\zeta(x)\dd x
 \eeqnn
for some constant $C_1>0$. As $G_1(y)\le G_2(y)$, we have $I^\de_{3,k}(s)\le 0$.
Combining the above inequalities with \eqref{5.6} and \eqref{5.5} we obtain
 \beqnn
 \ar\ar
\int_{\mbb{R}}\mbf{E}\big[\phi_k(v_{t\wedge\tau_n}^\delta(x))\big]J_\zeta(x)\dd x
-\int_{\mbb{R}}\phi_k(v_0^\delta(x))J_\zeta(x)\dd x \cr
 \ar\le\ar
\mbf{E}\Big[\int_0^{t\wedge\tau_n}\Big[\int_{\mbb{R}}
\big[2^{-1}C_1\zeta(\zeta+1)\phi_k(v_s^\delta(x)) \cr
 \ar\ar\qquad\qquad\qquad
+C\phi''_k(v_s^\delta(x))\<|v_s|,p_\delta^x\>^{2\gamma}\big]J_\zeta(x)\dd x
+I_{4,k}^\delta(s)\Big]\dd s\Big].
 \eeqnn
Letting $\delta\to0$ and using Lemma \ref{t2.5} and dominated convergence we get
 \beqnn
 \ar\ar
\int_{\mbb{R}}\mbf{E}\big[\phi_k(v_{t\wedge\tau_n}(x))\big]J_\zeta(x)\dd x
-\int_{\mbb{R}}\phi_k(v_0(x))J_\zeta(x)\dd x \cr
 \ar\le\ar
\mbf{E}\Big[\int_0^{t\wedge\tau_n}\dd s\int_{\mbb{R}}
\big[2^{-1}C_1\zeta(\zeta+1)\phi_k(v_s(x))
+C\phi''_k(v_s(x))|v_s(x)| \cdot|v_s(x)|^{2\gamma-1} \cr
 \ar\ar\qquad\qquad\qquad\qquad
+\phi'_k(v_s(x))[G_2(Y_s^{(1)}(x))-G_2(Y_s^{(2)}(x))]\big]J_\zeta(x)\dd x\Big].
 \eeqnn
Recalling the constant $r_0$ in Condition \ref{tc2}
and taking suitable $\zeta>0$ so that $r_0=2^{-1}C_1\zeta(\zeta+1)$.
Since $|v_0^+(x)|\equiv0$ and $0\le |y|\phi_k''(y)\le 2/k$ for all $y\in\mbb{R}$,  letting $k\to\infty$ in the above inequality we get
 \beqnn
 \ar\ar
\<\mbf{E}[v_{t\wedge\tau_n}^+],C_2 J_\zeta\>
=
C_2\int_{\mbb{R}}\mbf{E}\big[v_{t\wedge\tau_n}^+(x)\big]J_\zeta(x)\dd x \cr
 \ar\le\ar
C_2\mbf{E}\Big[\int_0^{t\wedge\tau_n}
\dd s\int_{\mbb{R}}[r_0v_s^+(x)
+G_2(Y_s^{(1)}(x))-G_2(Y_s^{(2)}(x))]1_{\{v_s(x)>0\}}J_\zeta(x)\dd x\Big] \cr
 \ar\le\ar
\int_0^t\dd s\int_{\mbb{R}}
\mbf{E}\big[r_n(v_{s\wedge\tau_n}^+(x))\big]C_2J_\zeta(x)\dd x
\le\int_0^tr_n(\<\mbf{E}[v_{s\wedge\tau_n}^+],C_2J_\zeta\>)\dd s,
 \eeqnn
where $C_2:=\<J_\zeta,1\>^{-1}$, $v_t^+(x):=v_t(x)\vee0$
and we used Condition \ref{tc2},
concaveness of $y\mapsto r_n(y)$ and Jensen's inequality in the last inequality.
It is then easy to show that
 \beqnn
\int_{\mbb{R}}\mbf{E}\big[v_{t\wedge\tau_n}^+(x)\big]J_\zeta(x)\dd x=0,
 \eeqnn
and hence, $\int_0^1v_t^+(x)\dd x=0$ for all $t\le \tau_n$ and $n\ge1$.
Letting $n\to\infty$ we obtain $\int_0^1v_t^+(x)\dd x=0$,
which ends the proof by the continuities of $(x,t)\mapsto Y_t^{(1)}(x)$
and $(x,t)\mapsto Y_t^{(2)}(x)$.
\qed

\section{Proof of Theorems \ref{t1.2} and \ref{t1.2'}}
\setcounter{equation}{0}

This section is devoted to the proofs of Theorems \ref{t1.2} and \ref{t1.2'}.
Recall that $\kappa_0:=\sup_{x,y\in[0,1]}\kappa(x,y)$.
Let $T>0$ and $p>1$ be fixed.
For simplicity we write $F_s^{t-s}(x):=F_s(p_{t-s}^x(\cdot)):=F_s(p_{t-s}(x,\cdot))$ for $F^{(i)}_s(p_{t-s}^{(i)}(x,\cdot))$ in the rest of this paper.

First we establish Theorem \ref{t1.2} under Lipschitz condition.

\bproposition\label{t3.5}
Suppose that Condition \ref{tc4.1}(i)--(ii) hold, and $G$ and $H$ satisfy the following
Lipschitz condition:
 \beqlb\label{5.10}
|G(x)-G(y)|+|H(x)-H(y)|\le C|x-y|,\qquad x,y\in\mbb{R}.
 \eeqlb
Then for each boundary conditions (C1) to (C4), \eqref{1.1} has a
strong solution $(Y_t)_{t\ge0}$.
\eproposition

We now proceed to proving Proposition \ref{t3.5}.
Let $(\tilde{Y}_t)_{t\ge0}$ be the solution to
\begin{equation}\label{3.2}
  \left\{
   \begin{aligned}
   \partial_t\tilde{Y}_t(x) &= \frac12\Delta \tilde{Y}_t(x),\quad t>0,~x\in(0,1),\\
   \tilde{Y}_0(x) &=Y_0(x),\qquad\qquad x\in[0,1]
   \end{aligned}
   \right.
  \end{equation}
with boundary condition (Ci), $i=1,2,3,4$.
For $n\ge0$, $t>0$ and $x\in[0,1]$
define $Y_t^{(n)}$  by the following:
 \beqlb\label{5.1}
Y_t^{(n)}(x)=\bar{Y}_t^{(n)}(x)+\tilde{Y}_t(x),
 \eeqlb
where $\bar{Y}_t^{(n)}$ are defined recursively by $\bar{Y}_t^{(0)}=0$ and
 \beqlb\label{3.1}
\bar{Y}_t^{(n+1)}(x)
 \ar=\ar
\int_0^t\dd s\int_0^1G(Y_s^{(n)}(y))p_{t-s}(x,y)\dd y \cr
 \ar\ar
+\int_0^t\int_0^1H(Y_s^{(n)}(y))p_{t-s}(x,y)W(\dd s,\dd y).
 \eeqlb
Recall that $P_tf(y)=\int_0^1p_t(x,y)f(x)\dd x$ for
$t>0$, $y\in[0,1]$ and $f\in\mathscr{B}[0,1]$.
For $n,t\ge0$ and $\delta>0$ define $\bar{Y}_t^{(n,\delta)}:=P_\delta \bar{Y}_t^{(n)}$.
It follows from \eqref{3.1} that
 \beqlb\label{3.1b}
\bar{Y}_t^{(n+1,\delta)}(x)
 \ar=\ar
\int_0^t\dd s\int_0^1G(Y_s^{(n)}(y))p_{t-s+\delta}(x,y)\dd y \cr
 \ar\ar
+\int_0^t\int_0^1H(Y_s^{(n)}(y))p_{t-s+\delta}(x,y)W(\dd s,\dd y).
 \eeqlb
Then for all $t,\delta>0$ and $n\ge1$,
$\bar{Y}_t^{(n,\delta)}\in C^2[0,1]$ and satisfies  boundary condition (Di).

\blemma\label{t4.5}
There is a constant $C_1>0$ so that
 \beqnn
\sup_{t\in[0,T]}\|\tilde{Y}_t\|_0\le C_1
\Big[\|Y_0\|_0+\sup_{t\in[0,T]}[|\mu_0(t)|+|\mu_1(t)|\Big].
 \eeqnn
\elemma
\proof
By Definition \ref{d1.1}, for each $f\in C^2[0,1]$ satisfying Condition (Di),
 \beqlb\label{4.2b}
\<\tilde{Y}_t,f\>
=
\<Y_0,f\>
+\frac12\int_0^t [\<\tilde{Y}_s,f''\>+F_s(f)]\dd s,
 \eeqlb
which can be written into the following mild form
 \beqnn
\<\tilde{Y}_t,f\>
=
\<Y_0,P_tf\>
+\frac12\int_0^t F_s(P_{t-s}f)\dd s.
 \eeqnn
It follows that
 \beqnn
\tilde{Y}_t(x)
=
\<Y_0,p_t^x\>
+\frac12\int_0^t F_s^{t-s}(x)\dd s.
 \eeqnn
Using Lemmas \ref{t2.2} and \ref{t5.2}  we obtain the assertion.
\qed

For $f\in\mathscr{B}[0,1]$ define $\|f\|^2:=\int_0^1|f(x)|^2\dd x$.
\blemma\label{t4.1}
Suppose that Condition \ref{tc4.1}(i)-(ii) holds.
Then
 \beqnn
\sup_{0\le t\le T,\,n\ge1}\mbf{E}\big[\|\bar{Y}_t^{(n)}\|^2+\|Y_t^{(n)}\|^2\big]<\infty.
 \eeqnn
\elemma
\proof
Applying H\"older's inequality, Lemma \ref{t2.5} and Condition \ref{tc4.1}(i), we get
 \beqlb\label{5.8}
 \ar\ar
\Big|\int_0^t\dd s\int_0^1G(Y_s^{(n)}(y))p_{t-s+\delta}(x,y)\dd y\Big|^2 \cr
 \ar\ar\quad\le
Ct\int_0^t\dd s\int_0^1|G(Y_s^{(n)}(y))|^2\cdot|p_{t-s+\delta}(x,y)|\dd y \cr
 \ar\ar\quad\le
Ct\int_0^t\dd s\int_0^1[|Y_s^{(n)}(y)|^2+1]\cdot|p_{t-s+\delta}(x,y)|\dd y
 \eeqlb
and
 \beqnn
M^{(n)}(t,s,x)
 \ar:=\ar
\Big|\int_0^1 \dd y\int_0^1H(Y_s^{(n)}(y))p_t(x,y)
H(Y_s^{(n)}(z))p_t(x,z)\kappa(y,z)\dd z\Big| \cr
 \ar\le\ar
\kappa_0 \Big[\int_0^1|H(Y_s^{(n)}(y))p_t(x,y)|\dd y\Big]^2
\le
C\int_0^1|H(Y_s^{(n)}(y))|^2\cdot|p_t(x,y)|\dd y \cr
 \ar\le\ar
C \int_0^1[|Y_s^{(n)}(y)|^2+1]\cdot|p_t(x,y)|\dd y.
 \eeqnn
It then follows from Doob's inequality that
 \beqnn
 \ar\ar
\mbf{E}\Big[\Big|\int_0^t\int_0^1H(Y_s^{(n)}(y))p_{t-s+\delta}(x,y)W(\dd s,\dd y)\Big|^2\Big]
\le
4\mbf{E}\Big[\int_0^tM^{(n)}(t-s+\delta,s,x)\dd s\Big] \cr
 \ar\ar\qquad\quad\le
C\mbf{E}\Big[\int_0^t\dd s\int_0^1[|Y_s^{(n)}(y)|^2+1]
\cdot|p_{t-s+\delta}(x,y)|\dd y\Big].
 \eeqnn
Now by using \eqref{3.1}, \eqref{5.8} and Lemma \ref{t4.5},
 \beqnn
\mbf{E}\big[\|Y_t^{(n+1,\delta)}\|^2\big]
\le
2\mbf{E}\big[\|\tilde{Y}_t\|^2\big]
+2\mbf{E}\big[\|\bar{Y}_t^{(n+1)}\|^2\big]
 \ar\le\ar
C+C\int_0^t\mbf{E}\big[\|Y_s^{(n)}\|^2+1\big]\dd s.
 \eeqnn
It follows from Lemma \ref{t2.5} and Fatou's lemma that
 \beqnn
\mbf{E}\big[\|Y_t^{(n+1)}\|^2\big]
 \ar=\ar
\mbf{E}\big[\|\lim_{\delta\to0}Y_t^{(n+1,\delta)}\|^2\big]
\le
\liminf_{\delta\to0}\mbf{E}\big[\|Y_t^{(n+1,\delta)}\|^2\big] \cr
 \ar\le\ar
C+C\int_0^t\mbf{E}\big[\|Y_s^{(n)}\|^2+1\big]\dd s.
 \eeqnn
Then by an induction argument one ends the proof.
\qed

\blemma\label{t4.2}
Under Condition \ref{tc4.1}(i)-(ii),
for each $f\in B[0,1]$, $\delta>0$ and $n\ge1$, we have
 \beqnn
\<\bar{Y}_t^{(n+1,\delta)},f\>
 \ar=\ar
\frac12\int_0^t
\<\Delta\bar{Y}_s^{(n+1,\delta)},f\>\dd s
+\int_0^t\<G(Y_s^{(n)}),P_\delta f\>\dd s \cr
 \ar\ar
+\int_0^t\int_0^1H(Y_s^{(n)}(y))P_\delta f(y)
W(\dd s,\dd y),
\qquad t\ge0.
 \eeqnn
\elemma
\proof
In view of \eqref{3.1b}, we deduce
 \beqnn
\<\bar{Y}_t^{(n+1,\delta)},f\>
=
\int_0^t\<G(Y_s^{(n)}),P_{t-s+\delta}f\>\dd s
+\int_0^t\int_0^1H(Y_s^{(n)}(y))P_{t-s+\delta}f(y)
W(\dd s,\dd y).
 \eeqnn
Put $k\ge1$ and $t_i:=it/k$. Then by virtue of Lemma \ref{t2.1}, for $i\ge1$,
 \beqlb\label{3.3}
 \ar\ar
\<\bar{Y}_{t_i}^{(n+1,\delta)},f\>
-\<\bar{Y}_{t_{i-1}}^{(n+1,\delta)},P_{t_i-t_{i-1}}f\> \cr
 \ar=\ar
\int_{t_{i-1}}^{t_i}\<G(Y_s^{(n)}),P_{t_i-s+\delta}f\>\dd s
+\int_{t_{i-1}}^{t_i}\int_0^1H(Y_s^{(n)}(y))P_{t_i-s+\delta}f(y)
W(\dd s,\dd y).
 \eeqlb
We deduce from \eqref{1.4} that
 \beqlb\label{4.17}
 \ar\ar
P_{t_i-t_{i-1}}f(y)-f(y)
=
\int_0^{t_i-t_{i-1}}\partial_s(P_sf(y))\dd s =
\frac12\int_0^{t_i-t_{i-1}}\dd s \int_0^1\Delta_x (p_s(x,y)) f(x)\dd x \cr
 \ar\ar\qquad=
\frac12\int_0^{t_i-t_{i-1}}\dd s \int_0^1\Delta_y (p_s(x,y)) f(x)\dd x
=
\frac12\int_0^{t_i-t_{i-1}}\Delta_y (P_sf(y))\dd s.
 \eeqlb
By the boundary conditions on $\bar{Y}_{t_{i-1}}^{(n+1,\delta)}$ and $p_{t_i-s}$,
 \beqnn
\tilde{G}_n(i,s):=\bar{Y}_{t_{i-1}}^{(n+1,\delta)}(y)\nabla_y (P_{t_i-s}f(y))\big|_0^1
-\nabla \bar{Y}_{t_{i-1}}^{(n+1,\delta)}(y) P_{t_i-s}f(y)\big|_0^1=0.
 \eeqnn
Applying integration by parts and \eqref{4.17} we get
 \beqnn
 \ar\ar
\<\bar{Y}_{t_{i-1}}^{(n+1,\delta)},P_{t_i-t_{i-1}}f-f\>
=\frac12\int_{t_{i-1}}^{t_i}
\<\bar{Y}_{t_{i-1}}^{(n+1,\delta)},\Delta(P_{t_i-s}f)\>\dd s \cr
 \ar\ar\quad=
\frac12\int_{t_{i-1}}^{t_i}
\Big[\<\Delta\bar{Y}_{t_{i-1}}^{(n+1,\delta)},P_{t_i-s}f\>
+\tilde{G}_n(i,s)\Big]\dd s
=
\frac12\int_{t_{i-1}}^{t_i}
\<\Delta\bar{Y}_{t_{i-1}}^{(n+1,\delta)},P_{t_i-s}f\>\dd s.
 \eeqnn
It thus follows from \eqref{3.3} that
 \beqnn
 \ar\ar
\<\bar{Y}_t^{(n+1,\delta)},f\> \cr
 \ar=\ar
\sum_{i=1}^k\big[\<\bar{Y}_{t_{i-1}}^{(n+1,\delta)},P_{t_i-t_{i-1}}f-f\>\big]
+
\sum_{i=1}^k\big[\<\bar{Y}_{t_i}^{(n+1,\delta)},f\>
-\<\bar{Y}_{t_{i-1}}^{(n+1,\delta)},P_{t_i-t_{i-1}}f\>\big] \cr
 \ar=\ar
\frac12\sum_{i=1}^k\int_{t_{i-1}}^{t_i}
\<\Delta\bar{Y}_{t_{i-1}}^{(n+1,\delta)},P_{t_i-s}f\>\dd s
+\sum_{i=1}^k\int_{t_{i-1}}^{t_i}
\<G(Y_s^{(n)}),P_{t_i-s+\delta}f\>\dd s  \cr
 \ar\ar
+\sum_{i=1}^k\int_{t_{i-1}}^{t_i}\int_0^1H(Y_s^{(n)}(y))P_{t_i-s+\delta}f(y)
W(\dd s,\dd y) \cr
 \ar=\ar
\frac12\int_0^t\sum_{i=1}^k 1_{(t_{i-1},t_i]}(s)
\<\Delta\bar{Y}_{t_{i-1}}^{(n+1,\delta)},P_{t_i-s}f\>\dd s \cr
 \ar\ar
+\int_0^t\sum_{i=1}^k 1_{(t_{i-1},t_i]}(s)
\<G(Y_s^{(n)}),P_{t_i-s+\delta}f\>\dd s \cr
 \ar\ar
+\int_0^t\int_0^1\sum_{i=1}^k 1_{(t_{i-1},t_i]}(s)H(Y_s^{(n)}(y))P_{t_i-s+\delta}f(y)
W(\dd s,\dd y).
 \eeqnn
Letting $k\to\infty$, and using Lemma \ref{t4.1}
and dominated convergence one can conclude the assertion.
\qed

\noindent{\it Proof of Proposition \ref{t3.5}}.
For $n,s\ge0$ and $\delta>0$ let
 \beqnn
u_s^{(n)}=\bar{Y}_s^{(n+1)}-\bar{Y}_s^{(n)}=Y_s^{(n+1)}-Y_s^{(n)},~~
u_s^{(n,\delta)}=P_\delta u_s^{(n)}
 \eeqnn
and
 \beqnn
H_s^{(n)}(y):=H(Y_s^{(n+1)}(y))-H(Y_s^{(n)}(y)),~
G_s^{(n)}(y):=G(Y_s^{(n+1)}(y))-G(Y_s^{(n)}(y)).
 \eeqnn
Then $u_s^{(n,\delta)}\in C^2[0,1]$ and satisfies boundary condition (Di) for each
fixed $s,n,\delta$.
For $\lambda,\delta>0$
let $v_s^{(n)}=\e^{-\lambda s}\|u_s^{(n)}\|^2$
and
$v_s^{(n,\delta)}=\e^{-\lambda s}\|u_s^{(n,\delta)}\|^2$.
In virtue of Lemma \ref{t4.2}, for each $f\in B^2[0,1]$.
 \beqnn
 \ar\ar
\<u_t^{(n,\delta)},f\>^2\e^{-\lambda t} \cr
 \ar=\ar
\int_0^t\<u_s^{(n,\delta)},f\>\<\Delta u_s^{(n,\delta)},f\>
\e^{-\lambda s}\dd s
+2\int_0^t\<u_s^{(n,\delta)},f\>\<G_s^{(n-1)},P_\delta f\>
\e^{-\lambda s}\dd s \cr
 \ar\ar
+\int_0^t\e^{-\lambda s}\dd s\int_0^1\dd y
\int_0^1H_s^{(n-1)}(y)P_\delta f(y)H_s^{(n-1)}(z)P_\delta f(z)
\kappa(y,z)\dd z \cr
 \ar\ar
+2\int_0^t\int_0^1\e^{-\lambda s}\<u_s^{(n,\delta)},f\>H_s^{(n-1)}(y)P_\delta f(y)
W(\dd s,\dd y)-\lambda\int_0^t\<u_s^{(n,\delta)},f\>^2\e^{-\lambda s}\dd s.
 \eeqnn
Summing on $f$ over a complete orthonormal system of $L^2[0,1]:=\{f\in\mathscr{B}[0,1]:\|f\|<\infty\}$ we obtain
 \beqlb\label{4.1b}
v_t^{(n,\delta)}
 \ar=\ar
\int_0^t\<u_s^{(n,\delta)},\Delta u_s^{(n,\delta)}\>
\e^{-\lambda s}\dd s
+2\int_0^t\<u_s^{(n,\delta)},P_\delta G_s^{(n-1)}\>\e^{-\lambda s}\dd s \cr
 \ar\ar
+\int_0^t\e^{-\lambda s}\dd s\int_0^1\dd y
\int_0^1H_s^{(n-1)}(y)H_s^{(n-1)}(z)p_{2\delta}(y,z)
\kappa(y,z)\dd z \cr
 \ar\ar
+2\int_0^t\int_0^1\e^{-\lambda s}H_s^{(n-1)}(y)u_s^{(n,2\delta)}(y)
W(\dd s,\dd y)-\lambda\int_0^tv_s^{(n,\delta)}\dd s.
 \eeqlb
Observe that under boundary conditions (D1) to (D4),
 \beqnn
u_s^{(n,\delta)}(0)\nabla u_s^{(n,\delta)}(0)
=u_s^{(n,\delta)}(1)\nabla u_s^{(n,\delta)}(1)=0.
 \eeqnn
Then by integration by parts,
 \beqnn
\<u_s^{(n,\delta)},\Delta u_s^{(n,\delta)}\>
=-\|\nabla u_s^{(n,\delta)}\|^2,
 \eeqnn
which deduces from \eqref{4.1b} that
 \beqlb\label{4.1c}
v_t^{(n,\delta)}
 \ar=\ar
-\int_0^t\|\nabla u_s^{(n,\delta)}\|^2
\e^{-\lambda s}\dd s
+2\int_0^t\<u_s^{(n,\delta)},P_\delta G_s^{(n-1)}\>\e^{-\lambda s}\dd s \cr
 \ar\ar
+\int_0^t\e^{-\lambda s}\dd s\int_0^1\dd y
\int_0^1H_s^{(n-1)}(y)H_s^{(n-1)}(z)p_{2\delta}(y,z)
\kappa(y,z)\dd z \cr
 \ar\ar
+2\int_0^t\int_0^1\e^{-\lambda s}H_s^{(n-1)}(y)u_s^{(n,2\delta)}(y)
W(\dd s,\dd y)-\lambda\int_0^tv_s^{(n,\delta)}\dd s.
 \eeqlb
It follows that
 \beqlb\label{5.4}
\mbf{E}[v_t^{(n,\delta)}]
 \ar\le\ar
2\int_0^t\e^{-\lambda s}\mbf{E}\big[\<u_s^{(n,\delta)},
P_\delta G_s^{(n-1)}\>\big]\dd s-\lambda\int_0^t\mbf{E}[v_s^{(n,\delta)}]\dd s
\cr
 \ar\ar
+\int_0^t\e^{-\lambda s}\dd s\int_0^1\dd y
\int_0^1\mbf{E}\big[H_s^{(n-1)}(y)H_s^{(n-1)}(z)\big]p_{2\delta}(y,z)\kappa(y,z)\dd z.
 \eeqlb
Under condition \eqref{5.10},
 \beqnn
|\<u_s^{(n,\delta)},P_\delta G_s^{(n-1)}\>|
\le
C\|u_s^{(n,\delta)}\|^2
+C\|u_s^{(n-1)}\|^2
 \eeqnn
and
 \beqnn
 \ar\ar
\int_0^1\dd y
\int_0^1\mbf{E}\big[H_s^{(n-1)}(y)H_s^{(n-1)}(z)\big]p_{2\delta}(y,z)\kappa(y,z)\dd z \cr
 \ar\ar\quad\le
2^{-1}\kappa_0
\int_0^1\dd y
\int_0^1\mbf{E}\big[|H_s^{(n-1)}(y)|^2+|H_s^{(n-1)}(z)|^2\big]p_{2\delta}(y,z)\dd z \cr
 \ar\ar\quad\le
C\int_0^1\mbf{E}\big[|H_s^{(n-1)}(y)|^2\big]\dd y
\le C\mbf{E}[\|u_s^{(n-1)}\|^2].
 \eeqnn
Combining above inequalities and \eqref{5.4} we obtain
 \beqnn
\mbf{E}\big[v_t^{(n,\delta)}\big]
\le
-\lambda\int_0^t \mbf{E}[ v_s^{(n,\delta)} ]\dd s
+C\int_0^t \mbf{E}[ v_s^{(n,\delta)}+v_s^{(n-1)} ]\dd s
\le
C\int_0^t \mbf{E}[v_s^{(n-1)} ]\dd s
 \eeqnn
for large $\lambda>0$.
From Fatou's lemma it follows that
 \beqnn
 \ar\ar
\mbf{E}\big[v_t^{(n)}\big]
\le
\liminf_{\delta\to0}\mbf{E}\big[v_t^{(n,\delta)}\big]
\le C\int_0^t \mbf{E}[ v_s^{(n-1)} ]\dd s.
 \eeqnn
Then it is elementary to see that
$(Y_t^{(n)})_{t\ge0}$ and
$(\bar{Y}_t^{(n)})_{t\ge0}$ are Cauchy
sequences, and $(Y_t)_{t\ge0}$ and $(\bar{Y}_t)_{t\ge0}$ denote the limits.
Moreover,
$Y_t(x)=\bar{Y}_t(x)+\tilde{Y}_t(x)$.
Letting $n\to\infty$ in \eqref{3.1} we can obtain
 \beqnn
\bar{Y}_t(x)
=
\int_0^t\dd s\int_0^1G(Y_s(y))p_{t-s}(x,y)\dd y
+\int_0^t\int_0^1H(Y_s(y))p_{t-s}(x,y)W(\dd s,\dd y).
 \eeqnn
It follows that
for any $f\in C^2[0,1]$
satisfying (Di), we have
 \beqnn
\<\bar{Y}_t,f\>
=\frac12\int_0^t\<\bar{Y}_s,f''\>\dd s
+\int_0^t\<G(Y_s),f\>\dd s
+\int_0^t\int_0^1H(Y_s(x))f(x)W(\dd s,\dd x),~
t\ge0.
 \eeqnn
In view of \eqref{4.2b}, one sees that $(Y_t)_{t\ge0}$
is a strong solution to \eqref{1.1}, which ends the proof.
\qed

Next we establish the weak existence part of Theorem \ref{t1.2} under a weaker condition.
\bproposition\label{t3.1}
Suppose that Condition \ref{tc4.1}(i)--(ii) holds.
Then for each boundary conditions (C1) to (C4), \eqref{1.1} has a weak solution $(Y_t)_{t\ge0}$.
\eproposition

We now proceed to proving Proposition \ref{t3.1}.
For $n\ge1$ and $x\in\mbb{R}$
let
 \beqnn
G_n(x)=\int_{\mbb{R}}q_{n^{-1}}(x-y)((G(y)\wedge n)\vee(-n))\dd y
 \eeqnn
and
 \beqnn
H_n(x)=\int_{\mbb{R}}q_{n^{-1}}(x-y)((H(y)\wedge n)\vee(-n))\dd y,
 \eeqnn
where the function $q_t(x)$ is defined in Section 2.
Then $\lim_{n\to\infty}G_n(x)=G(x)$ and
$\lim_{n\to\infty}H_n(x)= H(x)$ for all $x\in\mbb{R}$.
Moreover, for each fixed $n\ge1$,
$G_n$ and $H_n$ satisfy Lipschitz condition,
and hence by Proposition \ref{t3.5}, equation \eqref{1.1} with $G,H$
replaced by $G_n,H_n$ has a strong solution $(U_t^{(n)})_{t\ge0}$.
Let $(\tilde{Y}_t)_{t\ge0}$ be the solution to \eqref{3.2}
and
 \beqlb\label{5.2}
\bar{U}_t^{(n)}(x):=U_t^{(n)}(x)-\tilde{Y}_t(x),\qquad t\ge0,\,n\ge1,\, x\in\mbb{R}.
 \eeqlb
It follows that for $f\in C^2[0,1]$ and $t\ge0$,
 \beqlb\label{5.3b}
\<\bar{U}_t^{(n)},f\>
 \ar=\ar
\frac12\int_0^t\<\bar{U}_s^{(n)},f''\>\dd s
+\int_0^t\<G_n(U_s^{(n)}),f\>\dd s \cr
 \ar\ar
+\int_0^t\int_0^1H_n(U_s^{(n)}(x))f(x)W(\dd s,\dd x),
 \eeqlb
which implies that
 \beqlb\label{5.3}
\bar{U}_t^{(n)}(x)
 \ar=\ar
\int_0^t\dd s\int_0^1G_n(U_s^{(n)}(y))p_{t-s}(x,y)\dd y \cr
 \ar\ar
+\int_0^t\int_0^1H_n(U_s^{(n)}(y))p_{t-s}(x,y)W(\dd s,\dd y)
 \eeqlb
for $t\ge0$ and $x\in\mbb{R}$.

\blemma\label{t4.3}
Suppose that Condition \ref{tc4.1}(i)-(ii) holds.
Then
 \beqnn
\sup_{n\ge1}
\mbf{E}\Big[\sup_{0<t\le T}\|U_t^{(n)}\|^{2p}\Big]<\infty.
 \eeqnn
\elemma
\proof
Let $\bar{U}_t^{(n,\delta)}=P_\delta \bar{U}_t^{(n)}$
and $U_t^{(n,\delta)}=P_\delta U_t^{(n)}$.
Then $\bar{U}_t^{(n,\delta)}\in C^2[0,1]$ satisfies the boundary condition (Di).
As the same argument in \eqref{4.1c},
 \beqlb\label{4.1}
\|\bar{U}_t^{(n,\delta)}\|^2
 \ar=\ar
-\int_0^t\|\nabla\bar{U}_s^{(n,\delta)}\|^2\dd s
+2\int_0^t\<\bar{U}_s^{(n,\delta)},P_\delta G_n(Y_s^{(n)})\>\dd s \cr
 \ar\ar
+\int_0^t\dd s\int_0^1\dd y
\int_0^1H_n(U_s^{(n)}(y))H_n(U_s^{(n)}(z))p_{2\delta}(y,z)
\kappa(y,z)\dd z \cr
 \ar\ar
+2\int_0^t\int_0^1H_n(U_s^{(n)}(y))\bar{U}_s^{(n,2\delta)}(y)
W(\dd s,\dd y).
 \eeqlb
We finish the proof by the following three steps.

{\bf Step 1.}
For $\lambda,\delta>0,\,t\ge0$ and $n\ge1$ put
 \beqnn
 \ar\ar
f_n^\delta(t)=\|\nabla\bar{U}_t^{(n,\delta)}\|^2
\|\bar{U}_t^{(n,\delta)}\|^{2p-2},~k_n^\delta(t)=
\|\bar{U}_t^{(n,\delta)}\|^{2p-2}
\<\bar{U}_t^{(n,\delta)},P_\delta G_n(U_t^{(n)})\>,\cr
 \ar\ar
g_n^\delta(t)=\|\bar{U}_t^{(n,\delta)}\|^{2p-2}\int_0^1\dd y
\int_0^1H_n(U_t^{(n)}(y))H_n(U_t^{(n)}(z))p_{2\delta}(y,z)
\kappa(y,z)\dd z, \cr
 \ar\ar
h_n^\delta(t)=\|\bar{U}_t^{(n,\delta)}\|^{2p-4}
\int_0^1\dd y\int_0^1M_n^\delta(t,y)M_n^\delta(t,z)\kappa(y,z)\dd z,\cr
 \ar\ar
\cr
 \ar\ar
l_n^\delta(t)=
\int_0^t\int_0^1
\|\bar{U}_s^{(n,\delta)}\|^{2p-2}M_n^\delta(s,y)W(\dd s,\dd y),
 \eeqnn
where $M_n^\delta(t,y):=H_n(U_t^{(n)}(y))\bar{U}_t^{(n,2\delta)}(y)$.
In this step we estimate $g_n^\delta$, $h_n^\delta$ and $k_n^\delta$.

Under Condition \ref{tc4.1}(i),
 \beqnn
|G_n(x)|\le\int_{\mbb{R}}q_1(y)|G(x-y/\sqrt{n})|\dd y
\le\int_{\mbb{R}}q_1(y)(|x|+|y|/\sqrt{n}+1)\dd y
\le C(|x|+1)
 \eeqnn
and
 \beqlb\label{5.7}
|H_n(x)|\le C(|x|+1).
 \eeqlb
Then applying Lemmas \ref{t2.5} and \ref{t4.5}, Condition \ref{tc4.1}(ii), and \eqref{5.2}
one sees that
 \beqnn
2|\<\bar{U}_t^{(n,\delta)},P_\delta G_n(U_t^{(n)})\>|
 \ar\le\ar
\|\bar{U}_t^{(n,\delta)}\|^2+\|G_n(U_t^{(n)})\|^2 \cr
 \ar\le\ar
C[\|\bar{U}_t^{(n)}\|^2+\|U_t^{(n)}\|^2+1]
\le
C[\|\bar{U}_t^{(n)}\|^2+1]
 \eeqnn
and
 \beqnn
 \ar\ar
2\Big|\int_0^1\dd y
\int_0^1H_n(U_t^{(n)}(y))H_n(U_t^{(n)}(z))p_{2\delta}(y,z)
\kappa(y,z)\dd z\Big| \cr
 \ar\le\ar
\kappa_0\int_0^1\dd y
\int_0^1[H_n(U_t^{(n)}(z))^2+H_n(U_t^{(n)}(y))^2]|p_{2\delta}(y,z)|
\dd z
\le C[\|\bar{U}_t^{(n)}\|^2+1].
 \eeqnn
Similarly,
 \beqlb\label{4.6}
 \ar\ar
\Big|\int_0^1\dd y\int_0^1M_n^\delta(t,y)M_n^\delta(t,z)\kappa(y,z)\dd z\Big|^{1/2}
\le
\kappa_0^{1/2} \Big|\int_0^1|M_n^\delta(t,y)|\dd y\Big| \cr
 \ar\ar\quad\le
C\|H_n(U_t^{(n)})\|^2 +
C\int_0^1\dd y\int_0^1|\bar{U}_t^{(n)}(x)|^2\cdot|p_{2\delta}(x,y)|\dd x \cr
 \ar\ar\quad\le
C(\|H_n(U_t^{(n)})\|^2+ \|\bar{U}_t^{(n)}\|^2)
\le
C[\|\bar{U}_t^{(n)}\|^2+1].
 \eeqlb
Therefore,
 \beqlb\label{4.2}
g_n^\delta(t)
 \ar\le\ar
C(\|U_t^{(n)}\|^2+1)\|\bar{U}_t^{(n)}\|^{2p-2} \cr
  \ar\le\ar
C
\big[\|U_t^{(n)}\|^{2p}+\|\bar{U}_t^{(n)}\|^{2p}
+\|\bar{U}_t^{(n)}\|^{2p-2}\big]
\le
C
[\|\bar{U}_t^{(n)}\|^{2p}+1]
 \eeqlb
and
 \beqlb\label{4.3}
h_n^\delta(t)
\le
C(\|\bar{U}_t^{(n)}\|^2+1)^2
\|\bar{U}_t^{(n)}\|^{2p-4}
\le
C
[\|\bar{U}_t^{(n)}\|^{2p}+1],
 \eeqlb
where we used the fact $x^r\le x+1$ and
 \beqnn
xy\le p^{-1}x^p+(p-1)p^{-1}y^{p/(p-1)},\quad
(x+y)^{2p}\le 2px^{2p}+2py^{2p}
 \eeqnn
for all $x,y\ge0$ and $0\le r\le 1$.
Similarly,
 \beqlb\label{4.9}
k_n^\delta(t)
\le
C\|\bar{U}_t^{(n)}\|^{2p-2}\cdot(\|\bar{U}_t^{(n)}\|^2+1)
\le
C[\|\bar{U}_t^{(n)}\|^{2p}+1].
 \eeqlb

{\bf Step 2.}
We show that for all $p>1$,
 \beqlb\label{4.5}
\sup_{s\in[0,t],\,n\ge1}\mbf{E}\big[\|U_s^{(n)}\|^{2p}\big]<\infty,\qquad t>0.
 \eeqlb
It follows from \eqref{4.1} and It\^o's formula that
 \beqlb\label{4.4}
\|\bar{U}_t^{(n,\delta)}\|^{2p}
=
p\int_0^t\big[g_n^\delta(s)-f_n^\delta(s)+2(p-1)h_n^\delta(s)
+k_n^\delta(s)\big]\dd s
+pl_n^\delta(t).
 \eeqlb
Then by \eqref{4.2}-\eqref{4.9} and \eqref{4.4}
 \beqnn
\mbf{E}[\|\bar{U}_t^{(n,\delta)}\|^{2p}]
\le
C\int_0^t\mbf{E}\big[\|\bar{U}_s^{(n)}\|^{2p}+1\big]\dd s.
 \eeqnn
Using Lemma \ref{t2.5} and Fatou's lemma,
 \beqnn
\mbf{E}[\|\bar{U}_t^{(n)}\|^{2p}]
=\mbf{E}[\|\lim_{\delta\to0}\bar{U}_t^{(n,\delta)}\|^{2p}]
\le
\liminf_{\delta\to0}\mbf{E}[\|\bar{U}_t^{(n,\delta)}\|^{2p}]
\le
C\int_0^t\mbf{E}\big[\|\bar{U}_s^{(n)}\|^{2p}+1\big]\dd s.
 \eeqnn
Now by Gronwall's lemma we obtain \eqref{4.5}.

{\bf Step 3.}
By Doob's inequality and the arguments in \eqref{4.6} and \eqref{4.2},
for all $t\ge0$,
 \beqnn
 \ar\ar
\mbf{E}\Big[\sup_{s\in[0,t]}|l_n^\delta(s)|^2\Big]
 \le
4\mbf{E}\Big[\int_0^t\dd s\int_0^1\dd y
\int_0^1\|\bar{U}_s^{(n,\delta)}\|^{4p-4}|M_n^\delta(s,y)M_n^\delta(s,z)|
\kappa(y,z)\dd z\Big] \cr
 \ar\ar\qquad\qquad\le
C\mbf{E}\Big[\int_0^t\|\bar{U}_s^{(n)}\|^{4p-4}
(\|U_t^{(n)}\|^2+1)^2\dd s\Big]
\le
C\mbf{E}\Big[\int_0^t(\|\bar{U}_s^{(n)}\|^{4p}+1)\dd s\Big].
 \eeqnn
Combining the above inequality with \eqref{4.5} and \eqref{4.4} one can obtain
 \beqnn
\sup_{n\ge1,\,\delta>0}\mbf{E}\Big[\sup_{s\in[0,t]}
\|U_s^{(n,\delta)}\|^{2p}\Big]<\infty.
 \eeqnn
Now using Fatou's lemma,
 \beqnn
\sup_{n\ge1}\mbf{E}\Big[\sup_{s\in[0,t]}\|U_s^{(n)}\|^{2p}\Big]
 \ar\le\ar
\sup_{n\ge1}\mbf{E}\Big[\liminf_{\delta\to0}
\sup_{s\in[0,t]}\|U_s^{(n,\delta)}\|^{2p}\Big] \cr
 \ar\le\ar
\sup_{n\ge1}\liminf_{\delta\to0}\mbf{E}\Big[
\sup_{s\in[0,t]}\|U_s^{(n,\delta)}\|^{2p}\Big]
<\infty,
 \eeqnn
which finishes the proof.
\qed

\blemma\label{t4.4}
Suppose that Condition \ref{tc4.1}(i)-(ii) hold.
Then for each $\alpha\in(0,1)$,
 \beqnn
\mbf{E}\big[|\bar{U}_t^{(n)}(x_1)-\bar{U}_t^{(n)}(x_2)|^{2p}\big]
\le C|x_1-x_2|^{\alpha p},\qquad t\in(0,T],~x_1,x_2\in[0,1],~n\ge1
 \eeqnn
and
 \beqnn
\mbf{E}\big[|\bar{U}_{t+\varepsilon}^{(n)}(x)-\bar{U}_t^{(n)}(x)|^{2p}\big]
\le C\varepsilon^{\alpha p/2},\qquad\quad t\in(0,T],~\varepsilon>0,~x\in[0,1],~n\ge1.
 \eeqnn
\elemma
\proof
Since the proofs are similar, we only present the first one.
It follows from \eqref{5.3} that
 \beqlb\label{3.7}
 \ar\ar
\bar{U}_t^{(n)}(x_1)-\bar{U}_t^{(n)}(x_2) \cr
 \ar=\ar
\int_0^t\int_0^1M_n(x_1,x_2,t,s,y)W(\dd s,\dd y)
+\int_0^t\dd s\int_0^1G_n(U_s^{(n)}(y))k_{t-s}(x_1,x_2,y)\dd y \cr
 \ar=:\ar
I_t^n(x_1,x_2)+
J_t^n(x_1,x_2),
 \eeqlb
where
 \beqnn
k_t(x_1,x_2,y):=p_t(x_1,y)-p_t(x_2,y),~M_n(x_1,x_2,t,s,y)
:=
H_n(U_s^{(n)}(y))k_{t-s}(x_1,x_2,y)
.
 \eeqnn

By Lemma \ref{t2.2},
 \beqnn
|k_t(x_1,x_2,y)|
=
|k_t(x_1,x_2,y)|^{1-\alpha}|k_t(x_1,x_2,y)|^{\alpha}
\le Ct^{-(\alpha+1)/2}|x_1-x_2|^\alpha,
 \eeqnn
which deduces
 \beqlb\label{3.6}
\int_0^1|k_t(x_1,x_2,y))|^2\dd y
 \ar\le\ar
Ct^{-(\alpha+1)/2}|x_1-x_2|^\alpha
\int_0^1[|p_t(x_1,y)|+|p_t(x_2,y)|]\dd y \cr
 \ar\le\ar
Ct^{-(\alpha+1)/2}|x_1-x_2|^\alpha,
 \eeqlb
where the last inequality follows from Lemma \ref{t2.5}.
Applying H\"older's inequality, \eqref{5.7} and \eqref{3.6}, for all $n\ge1$,
 \beqnn
 \ar\ar
\Big|\int_0^1\dd y\int_0^1M_n(x_1,x_2,t,s,y)M_n(x_1,x_2,t,s,z)\kappa(y,z)\dd z\Big| \cr
 \ar\ar\qquad\le
\kappa_0\Big[\int_0^1|M_n(x_1,x_2,t,s,y)|\dd y\Big]^2 \cr
 \ar\ar\qquad\le
\kappa_0\|H_n(U_s^{(n)})\|^2\int_0^1|k_{t-s}(x_1,x_2,y)|^2\dd y \cr
 \ar\ar\qquad\le
C(\|U_s^{(n)}\|^2+1)
(t-s)^{-\frac{\alpha+1}{2}}|x_1-x_2|^\alpha.
 \eeqnn
Using Burkholder-Davis-Gundy's inequality we then get
 \beqnn
\mbf{E}\big[|I_t^n(x_1,x_2)|^{2p}\big]
 \ar\le\ar
C\mbf{E}\Big[\Big|\int_0^t\dd s\int_0^1\dd y\int_0^1M_n(x_1,x_2,t,s,y)M_n(x_1,x_2,t,s,z)\kappa(y,z)\dd z\Big|^p\Big] \cr
 \ar\le\ar
C|x_1-x_2|^{\alpha p}\mbf{E}\Big[\sup_{0\le s\le T}[\|U_s^{(n)}\|^{2p}+1]\Big]
\Big|\int_0^ts^{-(\alpha+1)/2}\dd s \Big|^p,~t\in(0,T].
 \eeqnn
Similarly,
 \beqnn
\mbf{E}\big[|J_t^n(x_1,x_2)|^{2p}\big]
\le
C|x_1-x_2|^{\alpha p}\mbf{E}\Big[\sup_{0\le s\le T}[\|U_s^{(n)}\|^{2p}+1]\Big]
\Big|\int_0^ts^{-(\alpha+1)/2}\dd s \Big|^p,~~t\in(0,T],
 \eeqnn
which yields the assertion by \eqref{3.7} and Lemma \ref{t4.3}.
\qed

\noindent{\it Proof of Proposition \ref{t3.1}}.
By Lemma \ref{t4.4} and Kolmogorov's criteria
(e.g. \cite[Corollary 16.9]{Kallenberg}), for each $T>0$,
the sequence
$\{\bar{U}_t^{(n)}(x):(t,x)\in[0,T]\times[0,1]\}$
is tight on $C([0,T]\times[0,1])$ and
hence, has a convergent subsequence $(\bar{U}_t^{(n_k)})_{t\ge0}$
which converges in law to $(\bar{U}_t)_{t\ge0}$.
It follows that $(U_t^{(n_k)})_{t\ge0}$
converges in law to $(U_t)_{t\ge0}$ as $k\to\infty$.
Thus
 \beqnn
M^{(n_k)}:=(U_t^{(n_k)},W_t,\mu_0(t),\mu_1(t))_{t\ge0}
\longrightarrow
M:=(U_t,W_t,\mu_0(t),\mu_1(t))_{t\ge0}
 \eeqnn
in law as $k\to\infty$.
Applying Skorokhod's representation, on another probability space,
there are continuous processes
 \beqnn
\hat{M}^{(n_k)}:=(\hat{U}_t^{(n_k)},
\hat{W}_t^{(n_k)},\hat{\mu}_0^{(n_k)}(t),\hat{\mu}_1^{(n_k)}(t))_{t\ge0}
~\mbox{ and }~
\hat{M}:=(\hat{U}_t,\hat{W}_t,\hat{\mu}_0(t),\hat{\mu}_1(t))_{t\ge0}
 \eeqnn
with the same distribution as $M^{(n_k)}$ and $M$, respectively.
Moreover,
$\hat{M}^{(n_k)}$ converges to
$\hat{M}$
almost surely.
In the following let $f\in C^2[0,1]$ satisfy (Di).
Since $(U_t^{(n)})_{t\ge0}$ is a strong solution to \eqref{1.1}
with $G_n$ and $H_n$ replaced by $G$ and $H$ by Proposition \ref{t3.5},
then
 \beqnn
\<U_t^{(n_k)},f\>
 \ar=\ar
\<Y_0,f\>
+
\frac12\int_0^t [\<U_s^{(n_k)},f''\>+F_s(f)]\dd s
+\int_0^t\<G_{n_k}(U_s^{(n_k)}),f\>\dd s \cr
 \ar\ar
+\int_0^t\int_0^1 H_{n_k}(U_s^{(n_k)}(y))f(y)W(\dd s,\dd y).
 \eeqnn
Define $\hat{F}_s^{(n)}(f)$ and $\hat{F}_s(f)$ similar as $F_s(f)$ with
$\mu_0,\mu_1$ replaced by $\hat{\mu}^{(n)}_0,\hat{\mu}^{(n)}_1$
and $\hat{\mu}_0,\hat{\mu}_1$, respectively.
Therefore,
 \beqnn
\<\hat{U}_t^{(n_k)},f\>
 \ar=\ar
\<Y_0,f\>
+\frac12\int_0^t [\<\hat{U}_s^{(n_k)},f''\>+\hat{F}_s^{(n_k)}(f)]\dd s
+\int_0^t\<G_{n_k}(\hat{U}_s^{(n_k)}),f\>\dd s \cr
 \ar\ar
+\int_0^t\int_0^1 H_{n_k}(\hat{U}_s^{(n_k)}(y))f(y)W^{n_k}(\dd s,\dd y).
 \eeqnn
It thus follows from \cite[Lemma 2.4]{XiongYang17}
that
 \beqnn
\<\hat{U}_t,f\>
 \ar=\ar
\<Y_0,f\>
+\frac12\int_0^t [\<\hat{U}_s,f''\>+ \hat{F}_s(f)]\dd s
+\int_0^t\<G(\hat{U}_s),f\>\dd s \cr
 \ar\ar
+\int_0^t\int_0^1 H(\hat{U}_s(y))f(y)W(\dd s,\dd y).
 \eeqnn
This ends the proof.
\qed

Now we are ready to proving

\noindent{\it Proof of Theorem \ref{t1.2}}.
The second assertion follows from Theorem \ref{t1.1} immediately.
Now let us prove the first one.
For each $n\ge1$ define stopping time
 \beqnn
\tau_n:=\inf\big\{t\ge0:\mu_0(t)+\mu_1(t)\ge n\big\}.
 \eeqnn
Then $\tau_n\to\infty$ almost surely as $n\to\infty$.
By Corollary \ref{t1.3} and Proposition \ref{t3.1} for each $\mu_1^n(t):=\mu_1(t\wedge\tau_n)$
and $\mu_0^n(t):=\mu_0(t\wedge\tau_n)$, there
is a unique strong solution $(Y_t^n)_{t\ge0}$ to \eqref{1.1} satisfying the boundary condition by \cite[Corollary 2.8]{Kurtz07}. Moreover, $Y_t^n=Y_t^{n+1}$ for all $t\le \tau_n$.
Let $Y_t=Y_t^n$ for each $t\le \tau_n$.
Then $(Y_t)_{t\ge0}$ satisfies \eqref{1.1} with the boundary condition
(C1) to (C4).
\qed

\noindent{\it Proof of Theorem \ref{t1.2'}}. The conclusions can be justified by using essentially
the same argument as those in the proofs of Theorems \ref{t1.1} and \ref{t1.2} with $p_t^{(1)}(x,y)$ and $p_t^{(2)}(x,y)$ replaced by
 \beqnn
p_t^{(1)}(x,y)=q_t(x,y)-q_t(-x,y),\quad
p_t^{(2)}(x,y)=q_t(x,y)+q_t(-x,y),
 \eeqnn
respectively. We omit the details.
\qed

\section{Proof of Theorem \ref{t6.1}}
\setcounter{equation}{0}
In this section we present the proof of Theorem \ref{t6.1}.
Let $(Y_t)_{t\ge0}$ be a weak solution to \eqref{1.1}
with the boundary condition (Ci) and $(p_t)_{t\ge0}$ be the solution to
\eqref{1.4} with the boundary condition (Di) for $i\in\{1,2,3,4\}$.
Write $F_s^{t-s}(x)$ for $F^{(i)}_s(p_{t-s}^{(i)}(x,\cdot))$.
Let arbitrary constants $T>0$ and $p\ge1$ be fixed and $Y_0\in C[0,1]$.
We will prove our theorem by a sequence of lemmas.

\blemma\label{t6.2}
For each $t>0$ and $f\in C[0,1]$,
 \beqlb\label{7.2}
\<Y_t,f\>
 \ar=\ar
\<Y_0,P_tf\>
+\int_0^t [2^{-1}\<F_s^{t-s},f\>+\<G(Y_s),P_{t-s} f\>]\dd s \cr
 \ar\ar
+\int_0^t\int_0^1H(Y_s(y))P_{t-s} f(y)W(\dd s,\dd y).
 \eeqlb
Moreover, for $x\in(0,1)$,
 \begin{eqnarray}\label{7.3}
Y_t(x)&=&
\<Y_0,p_t^x\>
+\int_0^t [2^{-1}F_s^{t-s}(x)+\<G(Y_s),p_{t-s}^x\>]\dd s\nonumber\\
&&+\int_0^t\int_0^1H(Y_s(y))p_{t-s}^x(y)W(\dd s,\dd y)
 \end{eqnarray}
with $p_t^x(y):=p_t(x,y)$.
\elemma
\proof
Let $t>0$ be fixed.
By Definition \ref{d1.1}, for each $f\in C[0,1]$,
 \beqlb\label{7.1}
\<Y_t,P_v f\>
 \ar=\ar
\<Y_0,P_v f\>+\frac12\int_0^t\dd s
\int_0^1\<Y_s,\Delta_x p_v (x,\cdot)\>f(x)\dd x
+\frac12\int_0^t\<F_s^v,f\>\dd s\cr
 \ar\ar
+\int_0^t\<G(Y_s),P_v f\>\dd s
+\int_0^t\int_0^1 H(Y_s(y))P_v f(y) W(\dd s,\dd y).
 \eeqlb
Set $t_i=it/n$ for $0\le i\le n$ and $n\ge1$.
In view of \eqref{1.4} and \eqref{7.1}, we obtain
 \beqnn
 \ar\ar
\<Y_t,f\>-\<Y_0,P_tf\> \cr
 \ar=\ar
\sum_{i=1}^n\<Y_{t_i},P_{t-t_i}f-P_{t-t_{i-1}}f\>
+\sum_{i=1}^n[\<Y_{t_i},P_{t-t_{i-1}}f\>-\<Y_{t_{i-1}},P_{t-t_{i-1}}f\>] \cr
 \ar=\ar
\frac12\sum_{i=1}^n\int_{t-t_{i-1}}^{t-t_i}\dd s
\int_0^1\<Y_{t_i},\partial_sp_s(x,\cdot)\>f(x)\dd x \cr
 \ar\ar
+\frac12\sum_{i=1}^n\int_{t_{i-1}}^{t_i}\dd s
\int_0^1\<Y_s,\Delta_x p_{t-t_{i-1}} (x,\cdot)\>f(x)\dd x \cr
 \ar\ar
+\sum_{i=1}^n\int_{t_{i-1}}^{t_i}
[2^{-1}\<F_s^{t-t_{i-1}},f\>+\<G(Y_s),P_{t-t_{i-1}} f\>]\dd s \cr
 \ar\ar
+\sum_{i=1}^n\int_{t_{i-1}}^{t_i}\int_0^1 H(Y_s(y))P_{t-t_{i-1}} f(y) W(\dd s,\dd y) \cr
 \ar=\ar
-\frac12\int_0^t\sum_{i=1}^n1_{(t_{i-1},t_i]}(s)
\int_0^1
\big[\<Y_{t_i},\Delta_xp_{t-s}(x,\cdot)\>
-\<Y_s,\Delta_x p_{t-t_{i-1}} (x,\cdot)\>\big]f(x)\dd x \cr
 \ar\ar
+\int_0^t\sum_{i=1}^n1_{(t_{i-1},t_i]}(s)
[2^{-1}\<F_s^{t-t_{i-1}},f\>+\<G(Y_s),P_{t-t_{i-1}} f\>]\dd s \cr
 \ar\ar
+\int_0^t\sum_{i=1}^n1_{(t_{i-1},t_i]}(s)
\int_0^1 H(Y_s(y))P_{t-t_{i-1}} f(y) W(\dd s,\dd y).
 \eeqnn
Letting $n\to\infty$ we obtain \eqref{7.2}.
Taking $f(y)=p_\delta^x(y)$ in \eqref{7.2} and then
letting $\delta\to\infty$ and using dominated convergence and
Lemma \ref{t2.5} we get \eqref{7.3}.
\qed

\blemma\label{t6.5}
Under Condition \ref{tc4.1}(i)-(ii), we have
 \beqnn
\sup_{0<t\le T,\,x\in[0,1]}\mbf{E}\big[|Y_t(x)|^{2p}\big]<\infty.
 \eeqnn
\elemma
\proof
For each $n\ge1$ define stopping time $\tau_n$ by
 \beqnn
\tau_n:=\inf\{t\ge0:\|Y_t\|_0\ge n\}.
 \eeqnn
Then $\tau_n\to\infty$ as $n\to\infty$.
Applying Lemma \ref{t6.2} we find
 \beqlb\label{7.9}
Y_t(x)1_{\{t\le \tau_n\}}
 \ar=\ar
\<Y_0,p_t^x\>1_{\{t\le \tau_n\}}
+2^{-1}1_{\{t\le \tau_n\}}\int_0^t F_s^{t-s}(x)\dd s+1_{\{t\le \tau_n\}}\int_0^t\<G(Y_s),p_{t-s}^x\>]\dd s \cr
 \ar\ar
+1_{\{t\le \tau_n\}}\int_0^t\int_0^1H(Y_s(y))p_{t-s}^x(y)W(\dd s,\dd y) \cr
 \ar=:\ar
I_t^{1,n}(x)+2^{-1}I_t^{2,n}(x)+I_t^{3,n}(x)+I_t^{4,n}(x).
 \eeqlb

By applying Lemma \ref{t2.5} we find
 \beqlb\label{7.10}
\mbf{E}\big[|I_t^{1,n}(x)|^{2p}\big]
\le
\|Y_0\|_0^{2p}\cdot\<1,|p_t^x|\>^{2p}
\le C\|Y_0\|_0^{2p}.
 \eeqlb
Using Lemmas \ref{t2.2} and \ref{t5.2} and Condition
\ref{tc4.1}(ii) one can conclude that
 \beqlb\label{7.11}
\sup_{n\ge1,\,0<t\le T,\,x\in[0,1]}\mbf{E}\big[|I_t^{2,n}(x)|^{2p}\big]
\le
\sup_{0<t\le T,\,x\in[0,1]}\mbf{E}\Big[\Big|\int_0^tF_s^{t-s}(x)\dd s\Big|^{2p}\Big]<\infty.
 \eeqlb
By Lemmas \ref{t2.5} and \ref{t2.2}, H\"older's inequality and Condition
\ref{tc4.1}(i),
 \beqlb\label{7.12}
\mbf{E}\big[|I_t^{3,n}(x)|^{2p}\big]
 \ar\le\ar
C\mbf{E}\Big[\int_0^t\<|G(Y_s)|^{2p},p_{t-s}^x\>]1_{\{s\le \tau_n\}}\dd s \Big] \cr
 \ar\le\ar
C\mbf{E}\Big[\int_0^t\<|Y_s|^{2p}1_{\{s\le \tau_n\}}+1,p_{t-s}^x\>]\dd s \Big] \cr
 \ar\le\ar
C\int_0^t\sup_{y\in[0,1]}\mbf{E}\big[|Y_s(y)1_{\{s\le \tau_n\}}|^{2p}\big]\dd s
+Ct.
 \eeqlb

Applying Lemmas \ref{t2.5} and \ref{t2.2}, H\"older's inequality and Condition
\ref{tc4.1}(i) we obtain
 \beqnn
M(t,s,x)
 \ar:=\ar
\Big|\int_0^1 \dd y \int_0^1H(Y_s(y))p_t^x(y)
H(Y_s(z))p_{t-s}^x(z)\kappa(y,z) \dd z\Big| \cr
 \ar\le\ar
\kappa_0\Big|\int_0^1 |H(Y_s(y))p_t^x(y)|\dd y\Big|^2
\le C\int_0^1 |H(Y_s(y))|^2\cdot|p_t^x(y)|\dd y \cr
 \ar\le\ar
C\int_0^1 [|Y_s(y)|^2+1]\cdot|p_t^x(y)|\dd y.
 \eeqnn
It then follows from H\"older's inequality again that
 \beqnn
\Big|\int_0^tM(t-s,s,x)1_{\{s\le \tau_n\}} \dd s\Big|^p
\le C\int_0^t\dd s\int_0^1 [|Y_s(y)1_{\{s\le \tau_n\}}|^{2p}+1]\cdot|p_{t-s}^x(y)|\dd y
 \eeqnn
for $0<t\le T$ and $n\ge1$.
Then by virtue of Burkholder-Davis-Gundy's inequality,
 \beqnn
\mbf{E}\big[|I_t^{4,n}(x)|^{2p}\big]
 \ar\le\ar
\mbf{E}\Big[\Big|1_{\{t\le \tau_n\}}\int_0^t\int_0^1H(Y_s(y))p_{t-s}^x(y)W(\dd s,\dd y)\Big|^{2p}\Big] \cr
 \ar=\ar
\mbf{E}\Big[\Big|1_{\{t\le \tau_n\}}\int_0^{t\wedge \tau_n}\int_0^1H(Y_s(y))p_{t-s}^x(y)W(\dd s,\dd y)\Big|^{2p}\Big] \cr
 \ar\le\ar
\mbf{E}\Big[\Big|\int_0^t\int_0^1H(Y_s(y))p_{t-s}^x(y)1_{\{s\le \tau_n\}}W(\dd s,\dd y)\Big|^{2p}\Big] \cr
 \ar\le\ar
C\mbf{E}\Big[\Big|\int_0^t M(t-s,s,x)1_{\{s\le \tau_n\}} \dd s\Big|^p\Big] \cr
 \ar\le\ar
C\int_0^t\sup_{x\in[0,1]}\mbf{E}\big[|Y_s(y)1_{\{s\le \tau_n\}}|^{2p}\big]\dd s
+C
 \eeqnn
for all $0<t\le T$ and $n\ge1$.
Combining this with \eqref{7.9} to \eqref{7.12} we have
 \beqnn
\sup_{x\in[0,1]}\mbf{E}\big[|Y_t(x)1_{\{t\le \tau_n\}}|^{2p}\big]
\le
C+C\int_0^t\sup_{x\in[0,1]}\mbf{E}\big[|Y_s(x)1_{\{s\le \tau_n\}}|^{2p}\big]\dd s,
\quad
0<t\le T,\,n\ge1.
 \eeqnn
Now applying Gronwall's lemma we know that
 \beqnn
\sup_{t\in(0,T],\,x\in[0,1],\,n\ge1}\mbf{E}\big[|Y_t(x)1_{\{t\le \tau_n\}}|^{2p}\big]
<\infty.
 \eeqnn
From Fatou's lemma it follows that
 \beqnn
\sup_{t\in(0,T],\,x\in[0,1]}\mbf{E}\big[|Y_t(x)|^{2p}\big]
 \ar=\ar
\sup_{t\in(0,T],\,x\in[0,1]}\mbf{E}\Big[\lim_{n\to\infty}|Y_t(x)1_{\{t\le \tau_n\}}|^{2p}\Big] \cr
 \ar\le\ar
\sup_{t\in(0,T],\,x\in[0,1],\,n\ge1}
\mbf{E}\Big[|Y_t(x)1_{\{t\le \tau_n\}}|^{2p}\Big]<\infty,
 \eeqnn
which completes the proof.
\qed

\blemma\label{t6.3}
Suppose that Condition \ref{tc4.1} holds.
Then
for boundary condition (C2),
 \beqnn
\mbf{E}\big[|Y_t(x_1)-Y_t(x_2)|^{2p}\big]
\le
C[t^{-p}+1]|x_1-x_2|^p,\qquad t\in(0,T],~
x_1,x_2\in[0,1]
 \eeqnn
and for boundary conditions (C1), (C3) and (C4),
 \beqnn
\mbf{E}\big[|Y_t(x_1)-Y_t(x_2)|^{2p}\big]
\le
C[t^{-p}+1]|x_1-x_2|^{2p\gamma_0},\quad t\in(0,T],~
x_1,x_2\in[0,1].
 \eeqnn
\elemma
\proof
Since the proofs are similar, we only state that of boundary condition (C3).
By applying Lemma \ref{t6.2},
 \beqlb\label{7.16}
 \ar\ar
Y_t(x_1)-Y_t(x_2) \cr
 \ar=\ar
\<Y_0,p_t^{x_1}-p_t^{x_2}\>
+\frac12\int_0^t [F_s^{t-s}(x_1)-F_s^{t-s}(x_2)]\dd s
+\int_0^t \<G(Y_s),p_{t-s}^{x_1}-p_{t-s}^{x_2}\>\dd s \cr
 \ar\ar
+\int_0^t\int_0^1H(Y_s(y))[p_{t-s}^{x_1}(y)-p_{t-s}^{x_2}(y)]W(\dd s,\dd y) \cr
 \ar=:\ar
\sum_{i=1,2,3,4}I_t^i(x_1,x_2),\qquad t\in(0,T],~
x_1,x_2\in[0,1].
 \eeqlb

In view of Lemmas \ref{t2.5} and \ref{t2.2},
 \beqnn
\<1,|p_t^{x_1}-p_t^{x_2}|\>
 \ar\le\ar
Ct^{-1/2}|x_1-x_2|^{1/2}\<1,\sqrt{|p_t^{x_1}-p_t^{x_2}|}\> \cr
 \ar\le\ar
Ct^{-1/2}|x_1-x_2|^{1/2}\sqrt{\<1,|p_t^{x_1}-p_t^{x_2}|\>}
\le
Ct^{-1/2}|x_1-x_2|^{1/2}
 \eeqnn
for all $t>0$, which leads to
 \beqlb\label{7.14}
|I_t^1(x_1,x_2)|^{2p}
\le
\|Y_0\|_0^{2p}\<1,|p_t^{x_1}-p_t^{x_2}|\>^{2p}
\le
Ct^{-p}|x_1-x_2|^p,\qquad t\in(0,T],
 \eeqlb
and by using Condition
\ref{tc4.1}(i) and H\"older's inequality we deduce
 \beqlb\label{7.17}
\mbf{E}\big[|I_t^3(x_1,x_2)|^{2p}\big]
 \ar\le\ar
C
\mbf{E}\Big[\Big|\int_0^t\<|Y_s|+1,|p_{t-s}^{x_1}-p_{t-s}^{x_2}|\>
\dd s\Big|^{2p}\Big] \cr
 \ar\le\ar
C\mbf{E}\Big[\int_0^t\<(|Y_s|+1)^{2p},|p_{t-s}^{x_1}-p_{t-s}^{x_2}|\>\dd s\Big]
\cdot\Big|\int_0^t\<1,|p_{t-s}^{x_1}-p_{t-s}^{x_2}|\>\dd s\Big|^{2p-1} \cr
 \ar\le\ar
Cc_0|x_1-x_2|^p,\qquad
t\in(0,T],
 \eeqlb
where
 \beqnn
c_0:=\sup_{0<t\le T,\,y\in[0,1]}\mbf{E}\big[|Y_s(y)|^{2p}+1\big]<\infty
 \eeqnn
by Lemma \ref{t6.5}.
Moreover, by Condition
\ref{tc4.1} (i),
 \beqnn
M_{t,s}(x_1,x_2)
 \ar:=\ar
\Big|\int_0^1\dd y\int_0^1H(Y_s(y))[p_t^{x_1}(y)-p_t^{x_2}(y)]H(Y_s(z))[p_t^{x_1}(z)-p_t^{x_2}(z)]
\kappa(y,z)\dd z\Big| \cr
 \ar\le\ar
C\Big|\int_0^1[|Y_s(y)|+1]|p_t^{x_1}(y)-p_t^{x_2}(y)|\dd y\Big|^2,
 \eeqnn
which leads to
 \beqnn
\mbf{E}\Big[\Big|\int_0^tM_{t-s,s}(x_1,x_2)\dd s\Big|^p\Big]
\le
C|x_1-x_2|^p,\qquad t\in(0,T].
 \eeqnn
Then by Burkholder-Davis-Gundy's inequality,
 \beqlb\label{7.18}
\mbf{E}\big[|I_t^4(x_1,x_2)|^{2p}\big]
\le
C\mbf{E}\Big[\Big|\int_0^tM_{t-s,s}(x_1,x_2)\dd s\Big|^p\Big]
\le
C|x_1-x_2|^p,\quad t\in(0,T].
 \eeqlb

Observe that
 \beqnn
F_s^{t-s}(x_1)-F_s^{t-s}(x_2)
=
\mu_0(s)[\nabla p_{t-s}^{x_1}(0)-\nabla p_{t-s}^{x_1}(0)]
+\mu_1(s)J_{t-s}(x_1,x_2)
 \eeqnn
with $J_t(x_1,x_2):=p_t^{x_1}(1)-p_t^{x_2}(1)$.
Then by using \eqref{2.4} and Lemma \ref{t5.3},
 \beqlb\label{7.15}
\mbf{E}\Big[\Big|\int_0^t\mu_0(s)[\nabla p_{t-s}^{x_1}(0)
-\nabla p_{t-s}^{x_1}(0)]\dd s\Big|^{2p}\Big]
\le
C[t^{-p}+1]|x_1-x_2|^{2p\gamma_0}.
 \eeqlb
In view of Lemma \ref{t2.2},
 \beqnn
|J_t(x_1,x_2)|
 \ar=\ar
|J_t(x_1,x_2)|^{1/2}\cdot|J_t(x_1,x_2)|^{1/2} \cr
 \ar\le\ar
Ct^{-1/4}t^{-1/2}|x_1-x_2|^{1/2}
=C t^{-3/4}|x_1-x_2|^{1/2},
 \eeqnn
which yields
 \beqnn
\mbf{E}\Big[\Big|\int_0^t|\mu_1(s)J_{t-s}(x_1,x_2)|\dd s\Big|^{2p}\Big]
\le
Ct^{p/2}\mbf{E}\Big[\sup_{s\in(0,T]}|\mu_1(s)|^{2p}\Big]|x_1-x_2|^p.
 \eeqnn
Combining the above inequality with \eqref{7.15} we obtain
 \beqnn
\mbf{E}\big[|I_t^2(x_1,x_2)|^{2p}\big]
\le
C[t^{-p}+1]|x_1-x_2|^{2p\gamma_0},\qquad t\in(0,T].
 \eeqnn
Together this with  \eqref{7.16}-\eqref{7.18} one ends the proof.
\qed

Similar argument as in the proof of Lemma \ref{t6.3},
we can deduce the following lemma and omit the proof.
\blemma\label{t6.4}
Suppose that Condition \ref{tc4.1} holds
and $0<T_1<T$.
Then for boundary condition (C2),
 \beqnn
\mbf{E}\big[|Y_{t_1}(x)-Y_{t_2}(x)|^{2p}\big]
\le
C|t_1-t_2|^{p/2},\qquad t_1,t_2\in[0,T],~
x\in[0,1]
 \eeqnn
and for boundary conditions (C1), (C3) and (C4),
 \beqnn
\mbf{E}\big[|Y_{t_1}(x)-Y_{t_2}(x)|^{2p}\big]
\le
C|t_1-t_2|^{p\gamma_0},\qquad t_1,t_2\in[T_1,T],~
x\in[0,1].
 \eeqnn
\elemma

\noindent{\it Proof of Theorem \ref{t6.1}}.
By Kolmogorov's continuity criteria (see e.g. \cite[Corollary 1.2(ii)]{Walsh}) and Lemmas \ref{t6.3} and \ref{t6.4},
the assertion follows immediately.
\qed


\begin{thebibliography}{99}

\bibitem{BallyGyPardoux94}
Bally, V., Gy\"ongy, I. and Pardoux, E. (1994):
White noise driven parabolic SPDEs with measurable
drift. \textit{J. Funct. Anal.} {\bf 120}, 484-510.

\bibitem{DalangMueller06}
Dalang, R. C., Mueller, C. and Zambotti, L. (2006):
Hitting properties of parabolic SPDEs with reflection.
\textit{Ann. Probab.} {\bf 34}, 1423-1450.

\bibitem{DalangZhang13}
Dalang, R. C. and Zhang, T. (2013):
H\"older continuity of solutions of SPDEs with reflection.
\textit{Commun. Math. Stat.} {\bf 1}, 133142.

\bibitem{D}
Dawson, D. A. (1993): Measure-valued Markov processes. \'Ecole
d'\'Et\'e de Probabilit\'es de Saint-Flour XXI---1991, 1--260, {\em
Lecture Notes in Math., \bf 1541}, Springer, Berlin.

\bibitem{Dawsonli06}
Dawson, D. A. and Li, Z. (2006): Skew convolution semigroups and affine Markov processes.
\textit{Ann. Probab.} \textbf{34}, 1103-1142.

\bibitem{Dawsonli12}
Dawson, D. A. and Li, Z. (2012): Stochastic equations, flows
and
measure-valued processes. \textit{Ann. Probab.} \textbf{40},
813-857.

\bibitem{Dawson}
Dawson, D. A., Vaillancourt, J. and Wang, H. (2000):
Stochastic partial differential equations for a class of interacting measure-valued diffusions.
\textit{Ann. Inst. Henri. Poincar\'e. Probab. Stat.} {\bf 36}, 167-180.

\bibitem{Donati-MartinPardoux93}
Donati-Martin, C. and Pardoux, E. (1993): White noise driven SPDEs with reflection. \textit{Probab. Theory Related Fields}. {\bf 95}, 1-24.

\bibitem{dynkin}
Dynkin, E. B. (1994): {\em An Introduction to Branching
Measure-Valued Processes. CRM Monograph Series, \bf 6}. American
Mathematical Society, Providence, RI.

\bibitem{Eth}
Etheridge, A. M. (2000): {\em An Introduction to Superprocesses}.
University Lecture Series 20. American Mathematical Society.

\bibitem{Gomez13}
Gomez, A., Lee, K., Mueller, C., Wei, A. and Xiong, J. (2013):
Strong uniqueness for an spde via backward doubly stochastic differential equations. \textit{Stat. Probab. Lett}. {\bf 83}, 2186-2190.

\bibitem{Gy95}
Gy\"ongy, I. (1995): On non-degenerate quasi-linear stochastic partial differential equations. \textit{Potential Anal.} {\bf 4}, 157-171.


\bibitem{GyPardoux93a}
Gy\"ongy, I. and Pardoux, E. (1993): On quasi-linear stochastic partial differential equations. \textit{Probab. Theory Related Fields}. {\bf 94}, 413-425.

\bibitem{GyPardoux93b}
Gy\"ongy, I. and Pardoux, E. (1993):
On the regularization effect of space-time white noise on quasi-linear parabolic partial differential equations.
\textit{Probab. Theory Related Fields}. {\bf 97}, 211-227.

\bibitem{HuLD}
Hu, Y., Lu, F. and Nualart, D. (2013): H\"older continuity of the solutions for a class of nonlinear SPDEs arising from one dimensional superprocesses.
\textit{Probab. Theory Relat. Fields.} {\bf 56}, 27-49.

\bibitem{Kalsi19}
Kalsi, J. (2019):
Existence of invariant measures for reflected stochastic
partial differential equations.
To appear in \textit{J. Theoret. Probab.}


\bibitem{Kallenberg}
Kallenberg, O. (2002): \textit{Foundations of modern probability}. Second edition.
Probability and its Applications (New York). Springer-Verlag, New York.

\bibitem{KoS88}
Konno, N. and Shiga, T. (1988): Stochastic
    partial
    differential equations for some measure-valued diffusions.
    \textit{Probab. Theory Related Fields}. \textbf{79},
    201-225.
\bibitem{Kurtz07}
Kurtz, T. G. (2007): The Yamada-Watanabe-Engelbert theorem for general stochastic equations and inequalities.
\textit{Electron. J. Probab}. \textbf{12}, 951-965.

\bibitem{Li11} Li, Z. (2011): \textit{Measure-Valued Branching Markov
    Process}. Springer, Berlin.

\bibitem{LiWXZ}
Li, Z., Wang, H., Xiong, J. and Zhou, X. (2012):
Joint continuity for the solutions to a class of nonlinear SPDEs.
\textit{Probab. Theory Relat. Fields.} {\bf 153}, 441-469.

\bibitem{Mitoma}
Mitoma, I. (1985): An $\infty$-dimensional inhomogeneous Langevin equation.
\textit{J. Funct. Anal.} \textbf{61}, 342-359.

\bibitem{MyP11} Mytnik, L. and Perkins, E. (2011): Pathwise
    uniqueness
    for stochastic heat equations with H\"older continuous
    coefficients:
    the white noise case. \textit{Probab. Theory Related
    Fields}. \textbf{149},
   1-96.

\bibitem{MPS06} Mytnik, L., Perkins, E. and Sturm, A. (2006):
On pathwise uniqueness for stochastic heat equations with
non-Lipschitz coefficients. \textit{Ann. Probab.} \textbf{34},
1910-1959.

\bibitem{N14} Neuman, E. (2018):
Pathwise uniqueness of the stochastic heat equations with
spatially inhomogeneous white noise. \textit{Ann. Probab.}
\textbf{46}, 3090-3187.


\bibitem{NualartPardoux92}
Nualart, D. and Pardoux, E. (1992): White noise driven quasilinear SPDEs with reflection. \textit{Probab. Theory Relat. Fields.} {\bf 93}, 77-89.

\bibitem{Per}
Perkins, E. A. (2002): {\em Dawson-Watanabe Superprocesses and
Measure-Valued Diffusions}. Ecole d'Et\'e de Probabiliy\'es de
Saint-Flour XXIX-1999. Lecture Notes Math. {\bf 1781}, 125-329.
Bernard, P. ed. Springer, Berlin.


\bibitem{Rei89} Reimers, M. (1989): One dimensional stochastic
    differential equations and the branching measure diffusion.
    \textit{Probab. Theory Related Fields}. \textbf{81},
    319-340.


\bibitem{RS13}
Rippl, T. and Sturm, A. (2013): New results on pathwise
uniqueness
for the heat equation with colored noise. \textit{Electron. J.
Probab.} \textbf{18}, 1-46.

\bibitem{Rosen87}
Rosen, J. (1987): Joint continuity of the intersection local times of Markov processes. \textit{Ann. Probab.} {\bf 15}, 659-675.

\bibitem{Sturm03}
Sturm, A. (2003): On convergence of population processes in
random environments to the stochastic heat equation with
colored
noise. \textit{Electron. J. Probab.} \textbf{8}, 267-272.

\bibitem{Walsh} Walsh, J. (1986): An introduction to stochastic partial differential equations.
\textit{Lecture Notes in Math}. \textbf{1180}, 266-439. Springer, Berlin.

\bibitem{X}
Xiong, J. (2013): {\em Three Classes of Nonlinear Stochastic Partial Differential Equations}. World Scientific.

\bibitem{Xio13} Xiong, J. (2013):
Super-Brownian motion as the unique strong solution to an SPDE.
\textit{Ann. Probab.} \textbf{41}, 1030-1054.


\bibitem{XiongYang17a}
 Xiong, J. and Yang, X. (2017):
 Strong existence and uniqueness to a class of nonlinear SPDEs driven by Gaussian colored noises. \textit{Stat. Probab. Lett}. {\bf 129}, 113-129.

\bibitem{XiongYang17} Xiong, J. and Yang, X. (2019): Existence and pathwise uniqueness to an SPDE driven by $\alpha$-stable colored noise.
    \textit{Stochastic Process. Appl.} {\bf 129}, 2681-2722.

\bibitem{XuZhang09}
Xu, T. and Zhang, T. (2009):
White noise driven SPDEs with reflection: existence, uniqueness and large
deviation principles. \textit{Stochastic Process. Appl.} {\bf 119}, 3453-3470.

\bibitem{YangZhang14}
Yang, J. and Zhang, T. (2014):
Existence and uniqueness of invariant measures for SPDEs with two reflecting
walls. \textit{J. Theoret. Probab.} {\bf 27}, 863-877.


\bibitem{Zhang10}
Zhang, T. (2010):
White noise driven SPDEs with reflection: Strong Feller properties and
Harnack inequalities. \textit{Potential Anal.} {\bf 33}, 137-151.

\bibitem{Zhang11}
Zhang, T. (2011): Systems of stochastic partial differential equations with reflection: existence and uniqueness.
\textit{Stochastic Process. Appl.} {\bf 121}, 1356-1372.










\end{thebibliography}
\end{document}